\documentclass{article}%
\usepackage{amsmath}
\usepackage{amssymb}
\usepackage{amsmath}
\usepackage{amsfonts}
\usepackage{graphicx}
\usepackage[latin1]{inputenc}
\usepackage[english]{babel}
\usepackage{dsfont}
\usepackage{hyperref}
\usepackage{url}
\usepackage{mathrsfs}
\usepackage{abstract} 
\usepackage{stmaryrd}
\usepackage{authblk}  
\DeclareMathOperator*{\esssup}{ess\,sup}

\usepackage{vmargin}
\setmarginsrb{2.5cm}{2.5cm}{2.5cm}{2.5cm}{0cm}{0cm}{1.25cm}{1.25cm}
\providecommand{\U}[1]{\protect\rule{.1in}{.1in}}
\newtheorem{theorem}{Theorem}

\newtheorem{lemma}[theorem]{Lemma}

\newtheorem{proposition}[theorem]{Proposition}
\newtheorem{remark}[theorem]{Remark}

\newenvironment{proof}[1][Proof]{\noindent\textbf{#1.} }{\ \rule{0.5em}{0.5em}}
\begin{document}
	\title{$\mathbb{L}^p$-solutions for Reflected BSDEs with jumps in a general filtration under stochastic Lipschitz coefficients}
	
%
	
	\author[1]{Yassine El Qalli\thanks{Email: \href{mailto:elqalli@insea.ac.ma}{elqalli@insea.ac.ma}}}

	\author[2]{Badr Elmansouri\thanks{Corresponding author. Emails: \href{mailto:badr.elmansouri@edu.uiz.ac.ma}{badr.elmansouri@edu.uiz.ac.ma} \, \& \, \href{mailto:b.elmansouri@uca.ac.ma}{b.elmansouri@uca.ac.ma}}\thanks{ORCID: \href{https://orcid.org/0000-0003-2603-894X}{ORCID: 0000-0003-2603-894X}}}

	\affil[1]{National Institute of Statistics and Applied Economics (INSEA), Rabat, B.P. 6217, Morocco.}
	
	\affil[2]{Cadi Ayyad University (UCA), National School of Applied Sciences of Marrakech (ENSA-M), BP 575, Avenue Abdelkrim Khattabi, 40000, Guéliz, Marrakech, Morocco.}

%
	
	%
%
	\date{}
	\maketitle
\begin{abstract}
In this paper, we establish existence and uniqueness of $\mathbb{L}^p$-solutions, for $p \in (1,2)$, to reflected backward stochastic differential equations (RBSDEs) in a general filtration supporting both a Brownian motion and an independent Poisson random measure. Our results are derived under suitable $\mathbb{L}^p$-integrability assumptions on the data and a stochastic Lipschitz condition on the coefficient.
\vspace{0.3cm}

\noindent \textbf{Keywords:} Reflected BSDEs, $\mathbb{L}^p$-solutions, general filtration, Stochastic Lipschitz coefficient.

\noindent \textbf{MSC 2020:} 60H10, 60H15, 60H30.
\end{abstract}

\section{Introduction}\label{intro}
Non-linear Backward Stochastic Differential Equations (BSDEs) were first introduced by Pardoux and Peng \cite{pardoux1990adapted}. Since then, the theory has developed considerably through its connections with mathematical finance \cite{el1997backward}, stochastic control and differential games \cite{EH}, and partial differential equations \cite{PP92}. A key line of research has been to relax the standard square-integrability assumptions on the terminal condition and the generator, leading to the study of $\mathbb{L}^p$ solutions for $p \in (1,2)$. The first analysis under a Lipschitz condition is due to El Karoui et al. \cite[Section~5]{el1997backward}, followed by Briand et al. \cite{BDHPS} who proved existence for $\mathbb{L}^p$ solutions with monotone generators and general growth in $y$. Subsequent works extended these results: Briand et al. \cite{BLS} for the one-dimensional case with linear growth in $z$; Chen \cite{Ch} for uniformly continuous generators; Ma et al. \cite{MFS} for monotone generators with general growth in $y$ and uniform continuity in $z$; and Tian et al. \cite{TJS} for one-dimensional BSDEs with discontinuous generators in $y$. Fan \cite{Fan2015} established results for multi-dimensional BSDEs under weak monotonicity, while Wang et al. \cite{wang2019lp} obtained minimal and maximal solutions under $p$-order weak monotonicity, general growth in $y$, and linear growth in $z$.  

The framework was subsequently extended to BSDEs with jumps by Tang and Li \cite{TL} and Rong \cite{rong1997solutions}, who included both a Brownian motion and a Poisson random measure. In more general filtrations, where the martingale representation property may fail \cite[Section~III.4]{jacod2003limit}, an additional orthogonal martingale term becomes necessary. For $\mathbb{L}^2$-solutions, this approach was first introduced by El Karoui and Huang \cite{ELH} and by Carbone et al. \cite{carbone2008backward} in the context of RCLL martingales. In the $\mathbb{L}^p$ setting with $p>1$, Kruse and Popier \cite{kruse2016bsdes,kruse2017lp} established existence and uniqueness for BSDEs driven by Brownian and Poisson components under monotonicity and integrability conditions. Yao \cite{Yao} studied the case $p \in (1,2)$ with non-Lipschitz generators, obtaining existence and uniqueness through a Lipschitz approximation method. El Jamali \cite{ElJamali2023} examined the case of a filtration generated by a Lévy process. More recently, Elmansouri and El Otmani \cite{ElmansouriOtmani2026} have considered generalized BSDEs with stochastic monotone generators in a filtration supporting a Brownian motion.

The notion of reflected BSDEs (RBSDEs) was introduced by El Karoui et al. \cite{EKPPQ}, where the solution is constrained to remain above a prescribed obstacle. Under square-integrability of both the terminal condition and the barrier, together with a Lipschitz generator, they proved existence and uniqueness in a Brownian filtration with continuous barriers. These assumptions were subsequently relaxed: Hamad\`{e}ne \cite{H} treated discontinuous barriers; Hamad\`{e}ne and Ouknine \cite{HMSOUK} extended the framework of \cite{EKPPQ} to Brownian motion combined with a Poisson random measure; and further results for RCLL obstacles were obtained in \cite{Ess,HO1}. For $\mathbb{L}^p$ solutions with $p \in (1,2)$, Hamad\`{e}ne and Popier \cite{BDHPS} established existence and uniqueness in a Brownian setting with Lipschitz generators, while Rozkosz and S{\l}omi{\'n}ski \cite{RS} obtained similar results under a monotonicity condition. More recently, Yao \cite{Yao1} studied RBSDEs with jumps and Lipschitz generators in $(y,z,v)$, proving existence and uniqueness of $\mathbb{L}^p$ solutions via a fixed-point method. Other contributions on $\mathbb{L}^p$ solutions of BSDEs with jumps include Klimsiak \cite{K,Kl}, Eddahbi et al. \cite{EFO}, and Elmansouri and Marzougue \cite{ELMANSOURI2025110407}.  

The aim of this note is to complement these works by studying RBSDEs with an RCLL obstacle having only totally inaccessible jumps, where the terminal value, the generator, and the obstacle are merely $p$-integrable ($p \in (1,2)$), and the driver satisfies a stochastic Lipschitz condition. This condition arises naturally in many financial applications, where the classical Lipschitz property may fail (see, e.g., \cite{elmansouri2025pricing,Elmansouri2025}).

The paper is organized as follows. Section~\ref{sec2} introduces the notations and assumptions, while Section~\ref{sec4} establishes existence and uniqueness of $\mathbb{L}^p$-solutions ($p \in (1,2)$) for reflected BSDEs via a penalization method.

\section{Notations, assumptions and preliminary results}\label{sec2}
Let $T>0$ be the time horizon and $(\Omega,\mathcal{F},(\mathcal{F}_{t})_{t\leq T},\mathbb{P})$ a filtered probability space whose filtration $(\mathcal{F}_{t})_{t\leq T}$ is complete, right-continuous, and quasi-left-continuous. We assume that $(\mathcal{F}_{t})_{t\leq T}$ \textit{supports} an $\mathbb{R}$-valued Brownian motion $(B_{t})_{t\leq T}$ and an independent martingale measure $\tilde{\mu}$ associated with a standard Poisson random measure $\mu$ on $\mathbb{R}^+ \times \mathcal{U}$, where $\mathcal{U} := \mathbb{R}^d \setminus \{0\}$ ($d>1$) is equipped with its Borel $\sigma$-algebra $\mathbb{U}$. The compensator of $\mu$ is $\nu(dt,de) = dt\,\lambda(de)$, with $\lambda$ a $\sigma$-finite measure on $\mathcal{U}$ satisfying $\int_{\mathcal{U}} (1 \wedge |e|^2)\,\lambda(de) < +\infty$, and such that, for every $\mathcal{G} \in \mathbb{U}$ with $\lambda(\mathcal{G}) < +\infty$, the process $\{\tilde{\mu}([0,t] \times \mathcal{G}) := (\mu - \nu)([0,t] \times \mathcal{G})\}_{t \leq T}$ is a martingale.

We will denote by
\begin{itemize}
	\item $[X, Y]$ (resp. $[X]$) the quadratic covariation (resp. quadratic variation) of given RCLL semimartingales $X$ and $Y$.
	\item $\mathcal{T}_{[t,T]}$ the set of stopping times $\tau$ such that $\tau \in [t,T]$.
	\item $\mathcal{P}$ the predictable $\sigma$-algebra on $\Omega \times [0,T]$.
\end{itemize}

Let $p>1$, $\beta > 0$ and $(a_{t})_{t\leq T}$ be a non-negative $\mathcal{F}_{t}$-adapted process. We define the increasing continuous process $A_t:= \int_0^t a^2_sds$, for all  $t\in  [0,T]$, and we introduce the following spaces:
\begin{itemize}
	\item[$\bullet$] $\mathbb{L}^p_\lambda$ is the set of $\mathbb{R}$-valued and $\mathbb{U}$-measurable mapping $\psi:\mathcal{U}\rightarrow \mathbb{R}$ such that
	$$\|\psi\|^p_{\mathbb{L}^p_\lambda}=\int_{\mathcal{U}}|\psi(e)|^p\lambda(de) < +\infty.$$
	\item[$\bullet$] $\mathcal{L}^{p}_\beta$ is the space of $\mathbb{R}$-valued and $\mathcal{F}_{T}$-measurable random variables $\xi$ such that
	$$\|\xi\|_{\mathcal{L}^{p}_\beta}=\left(\mathbb{E}\left[e^{\frac{p}{2}\beta A_T}|\xi|^{p}\right]\right)^\frac{1}{p} < +\infty.$$

	\item[$\bullet$] $\mathcal{S}^{p}_\beta$ is the space of $\mathbb{R}$-valued and $\mathbb{F}$-adapted RCLL processes $(Y_{t})_{t\leq T}$ such that\\
	$$\|Y\|_{\mathcal{S}^{p}_\beta}=\left(\mathbb{E}\left[\sup\limits_{0\leq t\leq T}e^{\frac{p}{2}\beta A_t}|Y_{t}|^{p}\right]\right)^\frac{1}{p} < +\infty.$$
	\item[$\bullet$] $\mathcal{S}^{p,A}_\beta$ is the space of $\mathbb{R}$-valued and $\mathbb{F}$-adapted RCLL processes $(Y_{t})_{t\leq T}$ such that\\
	$$\|Y\|_{\mathcal{S}^{p,A}_\beta}=\left(\mathbb{E}\left[\int_{0}^{T}e^{\frac{p}{2}\beta A_t}|Y_{t}|^{p}dA_t\right]\right)^\frac{1}{p} < +\infty.$$
	\item[$\bullet$] $\mathcal{H}^{p}_\beta$ is the space of $\mathbb{R}$-valued and $\mathbb{F}$-predictable processes $(Z_{t})_{t\leq T}$ such that\\
	$$\|Z\|_{\mathcal{H}^{p}_\beta}=\left(\mathbb{E}\left[\left(\int_{0}^{T}e^{\beta A_t}|Z_{t}|^{2}dt\right)^\frac{p}{2}\right]\right)^\frac{1}{p} < +\infty.$$
	\item[$\bullet$] $\mathfrak{L}^{p}_{\mu,\beta}$ is the space of $\mathbb{R}$-valued and $\mathcal{P}\otimes \mathbb{U}$-predictable processes $(V_t)_{t\leq T}$ such that\\
	$$\|V\|_{\mathfrak{L}^{p}_{\mu,\beta}}=\left(\mathbb{E}\left[\left(\int_0^T\int_{\mathcal{U}}e^{\beta A_t}|V_t(e)|^2\mu(dt,de)\right)^\frac{p}{2}\right]\right)^\frac{1}{p} < +\infty.$$
	\item[$\bullet$] $\mathcal{M}^p_\beta$ is the space of $\mathbb{R}$-valued RCLL martingales orthogonal to $B$ and $\mu$ such that\\
	$$\|M\|_{\mathcal{M}^{p}_{\beta}}=\left(\mathbb{E}\left[\left(\int_0^Te^{\beta A_t}d[M]_t\right)^\frac{p}{2}\right]\right)^\frac{1}{p} < +\infty.$$
	
	\item[$\bullet$] $\mathcal{S}^{p}$ is the space of $\mathbb{R}$-valued, continuous, increasing,  $\mathbb{F}$-adapted processes $(K_{t})_{t\leq T}$ such that
	$$\|K\|_{\mathcal{S}^{p}}=\left(\mathbb{E}\left[|K_{T}|^{p}\right]\right)^\frac{1}{p} < +\infty.$$
	
	\item[$\bullet$] $\mathfrak{B}^p_\beta:=\mathcal{S}^{p}_\beta\cap \mathcal{S}^{p,A}_\beta$ is a Banach space endowed with the norm
	$$\|Y\|_{\mathfrak{B}^{p}_\beta}^p=\|Y\|_{\mathcal{S}^{p}_\beta}^p+\|Y\|_{\mathcal{S}^{p,A}_\beta}^p.$$
	
	\item[$\bullet$] $\mathcal{E}^p_\beta=\mathfrak{B}^p_\beta \times \mathcal{H}^p_\beta \times \mathfrak{L}^{p}_{\mu,\beta} \times \mathcal{M}^p_\beta$.
\end{itemize}

For a given RCLL process $(\omega_t)_{t\leq T}$, $\omega_{t-}=\lim\limits_{s\nearrow t}\omega_s$, $t\leq T$ $(\omega_{0-}=\omega_0)$; $\omega_-:=(\omega_{t-})_{t\leq T}$ and $\Delta\omega_t=\omega_t-\omega_{t-}$.

\begin{remark}
	In what follows, the symbol $\mathsf{K}_\alpha$  denotes a generic positive $(\mathsf{K}_\alpha>0)$ constant that depend on a specific set of parameters $\alpha$ and may change from one line to another.
\end{remark}

We present a version of It\^{o}'s formula applied to the function $(t,x) \mapsto e^{\frac{p}{2}\beta A_t} |x|^p$ for $p \in (1,2)$, which is not sufficiently smooth. We set $\check{x} := |x|^{-1}x\mathds{1}_{\{x \neq 0\}}$.\\
This result, to be used repeatedly in what follows, is a slight modification of \cite[Lemma 7]{kruse2016bsdes}, where we add a predictable process of finite variation $(K_t)_{t \leq T}$. The case of a filtration generated by a Brownian motion is discussed in detail in \cite[Lemma 2.2]{BDHPS}. The proof is straightforward and is omitted, as it follows the same arguments as in\cite[Lemma 7]{kruse2016bsdes}.
\begin{lemma}\label{Ito}
	We consider the $\mathbb{R}$-valued semimartingale $(X_t)_{t\leq T}$ defined by
	$$ X_{t}=X_0+\int_{0}^{t}F_sds+\int_{0}^{t}Z_{s}dB_{s}+\int_0^t\int_{\mathcal{U}}V_s(e)\tilde{\mu}(ds,de)+M_t+K_t,$$
	such that:
	\begin{itemize}
		\item  $\mathbb{P}$-a.s. the process $K$ is predictable of bounded variation.
		\item $M$ is an RCLL local martingale orthogonal to both $B$ and $\mu$, i.e., $M \in \mathcal{M}_{\mathrm{loc}}$.
		\item $\left(F_t\right)_{t \leq T}$ is an $\mathbb{R}$-valued progressively measurable process and $\left(Z_t\right)_{t \leq T}$, $\left(V_t\right)_{t \leq T}$ are predictable processes with values in $\mathbb{R}$, such that $\int_0^T\left\{F_t+|Z_t|^2+\|V_t\|_\lambda^2\right\}dt<+\infty$, $\mathbb{P}$-a.s.
	\end{itemize}
	Then, for any $p\geq 1$ there exists a continuous and non-decreasing process $(\ell_t)_{t\leq T}$ such that
	\begin{equation*}
		\begin{split}\label{e0}
			&e^{\frac{p}{2}\beta A_t}|X_t|^p\nonumber\\
			&=|X_0|^p+\frac{p}{2}\beta\int_0^te^{\frac{p}{2}\beta A_s}|X_s|^pdA_s+\frac{1}{2}\int_0^te^{\frac{p}{2}\beta A_s}\mathds{1}_{\{p=1\}}d\ell_s+p\int_0^te^{\frac{p}{2}\beta A_s}|X_s|^{p-1}\hat{X}_sF_s ds\nonumber\\
			&+p\int_0^te^{\frac{p}{2}\beta A_s}|X_{s-}|^{p-1}\hat{X}_{s-}dK_s+p\int_0^te^{\frac{p}{2}\beta A_s}|X_s|^{p-1}\hat{X}_sZ_sdB_s\nonumber\\
			&+p\int_0^t\int_{\mathcal{U}}e^{\frac{p}{2}\beta A_s}|X_{s-}|^{p-1}\check{X}_{s-}V_s(e)\tilde{\mu}(ds,de)+c(p)\int_0^te^{\frac{p}{2}\beta A_s}|X_s|^{p-2}|Z_s|^2
			\mathds{1}_{\{X_s\neq0\}}ds\nonumber\\
			&+p\int_0^te^{\frac{p}{2}\beta A_s}|X_{s-}|^{p-1}\hat{X}_{s-}dM_s+c(p)\int_0^te^{\frac{p}{2}\beta A_s}|X_s|^{p-2}\mathds{1}_{\{X_s\neq0\}}d[M]^c_s\nonumber\\
			&+\int_0^t\int_{\mathcal{U}}e^{\frac{p}{2}\beta A_s}\left[|X_{s-}+V_s(e)|^{p}-|X_{s-}|^p-p|X_{s-}|^{p-1}\check{X}_{s-}V_s(e)\right]\mu(ds,de)\nonumber\\
			&+\sum_{0 < s \leq t} e^{\frac{p}{2}\beta A_s}\left[|X_{s-}+ \Delta M_s|^{p}-|X_{s-}|^p-p|X_{s-}|^{p-1}\check{X}_{s-}\Delta M_s\right],
			\nonumber\\
		\end{split}
	\end{equation*}
	where $c(p)=\frac{p(p-1)}{2}$ and  $(\ell_t)_{t\leq T}$ is a continuous, non-decreasing process that increases only on the boundary of the random set $\{t\in  [0,T],\;\; X_{t-}=X_t=0\}$.
\end{lemma}

In  this paper, we aim to study reflected BSDE of the following form:
\begin{equation}\label{basic RBSDE}
	\displaystyle\left\{
	\begin{split}
		Y_{t}&=\xi +\displaystyle\int_{t}^{T} f(s,Y_{s},Z_{s},V_s)ds+(K_{T}-K_{t})-\displaystyle\int_{t}^{T}Z_{s}dB_{s}-\displaystyle\int_t^T\displaystyle\int_{\mathcal{U}}V_s(e)\tilde{\mu}(ds,de)\\
		&\quad\quad-\int_{t}^{T}dM_s,~t \in [0,T], \hbox{}\\
		Y_t &\geq L_t \quad \forall t \in [0,T],~~ \text{ and } \displaystyle\int_0^T(Y_t-L_t)dK_t=0\; \text{ a.s.}~~ \hbox{  }
	\end{split}
	\right.
\end{equation}
The problem consists of finding a quintuplet of $\mathbb{F}$-adapted processes $(Y,Z,V,M,K) \in \mathcal{E}^p_\beta \times \mathcal{S}^p$, for $p \in (1,2)$, that satisfies \eqref{basic RBSDE}.

\paragraph*{Assumptions on the data $(\xi,f,L)$}
We consider that
\begin{description}
	\item[$(\mathcal{H}1)$] The terminal condition $\xi \in \mathcal{L}^{p}_\beta$.
	\item[$(\mathcal{H}2)$] The coefficient $f : \Omega\times [0,T]\times\mathbb{R}\times\mathbb{R}\times (\mathbb{L}^1_\lambda+  \mathbb{L}^2_\lambda) \longrightarrow \mathbb{R}$ satisfies :
	\begin{itemize}
		\item[(i)] for all $(t,y,z,v)\in [0,T]\times\mathbb{R}^{}\times\mathbb{R}^{}\times(\mathbb{L}^1_\lambda+  \mathbb{L}^2_\lambda)$, $(\omega,t) \mapsto f(\omega,t,y,z,v)$ is an $\mathbb{F}$-progressively measurable process.
		\item[(ii)] There exists three positive $\mathbb{F}$-progressively measurable processes $(\theta_t)_{t\leq T}$, $(\gamma_t)_{t\leq T}$ and $(\eta_t)_{t\leq T}$ such that for all $(t,y,y',z,z',v,v')\in [0,T]\times\mathbb{R}^{2}\times\mathbb{R}^{2}\times(\mathbb{L}^1_\lambda+  \mathbb{L}^2_\lambda)^2$ ,
		$$|f(t,y,z,v)-f(t,y',z',v')|\leq   \theta_t |y-y'|+\gamma_t |z-z'|+\eta_t\|v-v'\|_{\mathbb{L}^1_\lambda+  \mathbb{L}^2_\lambda}.$$
		\item[(iii)] There exists $\epsilon>0$ such that $a^2_t:=\theta_t+\gamma^q_t+\eta^q_t \geq \epsilon$ with $q=\frac{p}{p-1}$ for each $p\in(1,2)$.
		\item[(iv)] The $\mathcal{F}_T$-measurable random variable $A_T$ is bounded by some constant $\mathfrak{C}$.
		\item[(v)] For all $(y,z,v)\in \mathbb{R} \times \mathbb{R}^{d}\times\mathfrak{L}_\lambda$, the process $(f(t,y,z,v))_{t\leq T}$ is progressively measurable and 
		$$\mathbb{E}\left[\left(\int_{0}^{T}e^{\beta A_t}\left|\frac{f(t,0,0,0)}{a_t}\right|^{2}dt\right)^\frac{p}{2}\right]<+\infty.$$
	\end{itemize}
	
	\item[$(\mathcal{H}3)$] The obstacle $(L_{t})_{t\leq T}$ is an RCLL, progressively measurable, real-valued process satisfying:
	\begin{itemize}
		\item[\textnormal{(i)}] $L_T\leq \xi$.
		\item[\textnormal{(ii)}] $\mathbb{E}\left[\sup\limits_{0\leq t\leq T}\left|e^{\frac{q}{2}\beta A_t}L_t^+\right|^{p}\right]<+\infty$, where $L^+$ denotes the positive part of $L$.
		\item[\textnormal{(iii)}] The jump times of $L$ are assumed to be inaccessible stopping times.
	\end{itemize}
\end{description}
\begin{remark}
	\begin{itemize}
		\item The boundedness condition $(\mathcal{H}2)$-\textnormal{(iii)} is imposed purely for technical reasons: it will be used in the proofs to derive the fundamental estimates and to ensure that the process $A_t=\int_0^t a_s^2\,\mathrm{d}s$ remains well behaved on $[0,T]$, in particular guaranteeing that the $q$-th powers of $\gamma$ and $\eta$ (with $q=\tfrac{p}{p-1}$) stay uniformly controlled as $p$ approaches $1$. Note that this assumption does \emph{not} imply that $\theta$, $\gamma$ or $\eta$ are themselves bounded: for example, if $\theta_t=1/\sqrt{t}$ on $(0,T]$, then $\int_0^T\theta_t\,\mathrm{d}t=2\sqrt{T}<+\infty$ while $\sup_{t\in(0,T]}\theta_t=+\infty$. Thus, assumptions $(\mathcal{H}2)$-\textnormal{(i)} and $(\mathcal{H}2)$-\textnormal{(iii)} do not imply that the driver $f$ is Lipschitz in the usual (deterministic) sense, but only that it is stochastically Lipschitz.
		
		\item As mentioned in \cite[page 4]{HMSOUK}, assumption $(\mathcal{H}3)$-\textnormal{(iii)} is satisfied if, for example, $\forall t \in [0,T]$, $L_t(\omega)=\mathcal{L}_t(\omega)+\mu(\omega,t,\mathfrak{S})$, where $\mathcal{L}$ is
		continuous and $\mathfrak{S}$ is a Borel set such that $\lambda(\mathfrak{S})<+\infty$.
	\end{itemize}
\end{remark}

\subsubsection*{Comparison theorem}
The comparison principal will be established under a slightly modified assumption on the driver $f$ with respect to the jump parameter $v$ in order to obtain a comparison result. More precisely, we retain assumptions $(\mathcal{H}1)$ and $(\mathcal{H}2)$, but replace condition $(\mathcal{H}2)$-\textnormal{(ii)} on the jump parameter $v$ with a monotonicity assumption described as follows:
\begin{description}
	\item[$(\mathcal{H}2)$]
	\begin{itemize}
		\item[(ii')] 
		\begin{itemize}
			\item[(a)] There exists three strictly positive $\mathcal{F}_t$-adapted processes $(\theta_t)_{t\leq T}$, $(\gamma_t)_{t\leq T}$ and $(\eta_t)_{t\leq T}$ such that for all $(t,y,y',z,z',v)\in [0,T]\times\mathbb{R}^{2}\times\mathbb{R}^{2}\times \mathbb{L}^1_\lambda+  \mathbb{L}^2_\lambda$
			$$|f(t,y,z,v)-f(t,y',z',v)|\leq   \theta_t |y-y'|+\gamma_t |z-z'|.$$
			
			\item[(b)]  For each $(y, z, \psi, \phi) \in \mathbb{R} \times$ $\mathbb{R} \times( \mathbb{L}^1_\lambda+  \mathbb{L}^2_\lambda) \times (\mathbb{L}^1_\lambda+  \mathbb{L}^2_\lambda)$, there exists a predictable process $\kappa=\kappa^{y, z, \psi, \phi}: \Omega \times[0, T] \times \mathcal{U} \rightarrow \mathbb{R}$ such that:
			$$
			f(t, y, z, \psi)-f(t, y, z, \phi) \leq \int_{\mathcal{U}}(\psi(e)-\phi(e)) \kappa_t^{y, z, \psi, \phi}(e) \pi({d} e)
			$$
			with $\mathbb{P} \otimes dt \otimes \lambda$-a.e. for any $\left(y, z, \psi, \phi\right)$,
			\begin{itemize}
				\item $\kappa_t^{y, z, \psi, \phi}(e) \geq -1$
				
				\item $\left|\kappa_t^{y, z, \psi, \phi}(e)\right| \leq \vartheta(e)$ where $\vartheta \in \mathbb{L}_\lambda^\infty \cap \mathbb{L}_\lambda^2$ and $\|\vartheta\|_{\mathbb{L}_\lambda^\infty \cap \mathbb{L}_\lambda^2} \leq \eta_t$ for all $t \in [0,T]$.
			\end{itemize}
		\end{itemize}
	\end{itemize}
\end{description}

Note that, since $\mathbb{L}_\lambda^\infty \cap \mathbb{L}_\lambda^2$ is the dual space of $\mathbb{L}^1_\lambda+  \mathbb{L}^2_\lambda$ (see Theorem 3.1 in \cite[Ch II]{KreinPetuninSemenov1982}), assumption $(\mathcal{H}2)$-\textnormal{(ii')}-\textnormal{(b)} yields,
for all $(t,y,z,v,v')\in [0,T]\times\mathbb{R}\times\mathbb{R}\times(\mathbb{L}^1_\lambda+  \mathbb{L}^2_\lambda)^2$,
$$
|f(t,y,z,v)-f(t,y,z,v')|\leq \|\vartheta\|_{\mathbb{L}_\lambda^\infty \cap \mathbb{L}_\lambda^2} \|v-v'\|_{\mathbb{L}^1_\lambda+  \mathbb{L}^2_\lambda} \leq    \eta_t\|v-v'\|_{\mathbb{L}^1_\lambda+  \mathbb{L}^2_\lambda}.
$$
Thus, under this assumption, we recover the classical stochastic Lipschitz property of $f$ given by $(\mathcal{H}2)$-\textnormal{(ii)}. Moreover, using this monotonicity assumption on $v$ together with the boundedness of $A_T$, we can derive the following comparison principle stated in \cite[Proposition 4]{kruse2017lp}.
\begin{proposition}\label{com}
	Let $f_1$ and $f_2$ be generators satisfying $(\mathcal{H}2')$. Let $\xi^1$ and $\xi^2$ be two terminal conditions for BSDEs (2) driven respectively by $f_1$ and $f_2$. Denote by $(Y^1, Z^1, V^1, M^1)$ and $(Y^2, Z^2, V^2, M^2)$ their respective solutions in some space $\mathcal{E}^p(0, T)$ with $p>1$. If $\xi^1 \leq \xi^2$ and $f_1\left(t, Y_t^1, Z_t^1, V_t^1\right) \leq f_2\left(t, Y_t^1, Z_t^1, V_t^1\right)$, then a.s., for any $t \in[0, T]$, $Y_t^1 \leq Y_t^2$.
\end{proposition}

\section{Existence and uniqueness of an $\mathbb{L}^p$-solutions}\label{sec4}
The main result of this paper is stated as follows:
\begin{theorem}\label{thm2}
	Under assumptions $(\mathcal{H}1)$, $(\mathcal{H}2)$ and $(\mathcal{H}3)$, the reflected BSDE \eqref{basic RBSDE} admits a unique $\mathbb{L}^p$-solution $(Y,Z,V,M,K) \in \mathcal{E}^p_\beta \times \mathcal{S}^p$ for $p \in (1,2)$. 
\end{theorem}

The proof of Theorem \ref{thm2} is divided into two main steps:
\begin{enumerate}
	\item \textbf{Step 1:} We first treat the case where $f$ does not depend on the solution, using a penalization method.
	
	\item \textbf{Step 2:} We then extend the result to the general stochastic Lipschitz case via a fixed-point argument in a suitable Banach space.
\end{enumerate}
The existence and uniqueness argument in the classical Brownian setting can be found in \cite[Section 4]{BDHPS}.

\begin{proof}
	\emph{}
	\paragraph*{Part 1: Case of a driver $f$ independent of the parameters $(y,z,v)$}\emph{}\\
	In this part, we aim to prove the existence result when $f$ does not depend on the jump variable $v$, meaning that
	$$
	f(\omega,t,y,z,v)=f(\omega,t,0,0,0)=:\mathfrak{f}(\omega,t),
	$$
	for all $(\omega,t,y,z,v)\in \Omega \times [0,T]\times\mathbb{R}\times\mathbb{R}\times(\mathbb{L}^1_\lambda+  \mathbb{L}^2_\lambda)$. We then establish existence and uniqueness for the following RBSDE:
	\begin{eqnarray}\label{basic RBSDE thm}
		&&\hspace{-1cm} \left\{
		\begin{array}{ll}
			Y_{t}=\xi +\displaystyle\int_{t}^{T} \mathfrak{f}(s)ds+K_{T}-K_{t}-\displaystyle\int_{t}^{T}Z_{s}dB_{s}-\displaystyle\int_t^T\displaystyle\int_{\mathcal{U}}V_s(e)\tilde{\mu}(ds,de)-\int_{t}^{T}dM_s,~t \in [0,T],& \hbox{}\\[0.3cm]
			Y_t\geq L_t, \quad \forall t \in [0,T],~\text{ and } \displaystyle\int_0^T(Y_t-L_t)dK_t=0 \quad \text{ a.s.} & \hbox{}
		\end{array}
		\right.
	\end{eqnarray}
	
	The proof follows a penalization approach. To implement it, we adapt the argumentation of Theorem 2 in \cite{kruse2017lp} (see also \cite[Theorem 8]{elmansouri2025mathbblpsolutionsbsdesreflectedbsdes} for a detailed exposition, which itself is inspired by the reasoning of Proposition 3 and Theorem 3 in \cite{kruse2017lp}). For each $n \in \mathbb{N}$, there exists a unique quadruple $(Y^n,Z^n,V^{n},M^n) \in \mathcal{E}^{p}_\beta$ (with $p \in (1,2)$) satisfying the following penalized BSDEs:
	\begin{equation}\label{penalized}
		Y^n_{t}=\xi +\int_{t}^{T} \mathfrak{f}(s)ds+n\int_{t}^{T}(Y^n_s-L_s)^-ds-\int_{t}^{T}Z^n_{s}dB_{s}
		-\int_t^T\int_{\mathcal{U}}V_s^n(e)\tilde{\mu}(ds,de)-\int_{t}^{T} dM^n_s.
	\end{equation}
	Define $K^n_t := n\int_{0}^{t}(Y^n_s - L_s)^- \, ds$. The proof of \textbf{Part 1} consists of four steps.
	
	~\\
	\textbf{Step 1:} There exists a positive constant $\mathsf{K}_{p,\mathfrak{C},\beta}>0$ (independent on $n$) such that
	\begin{eqnarray*}
		&&\hspace{-1cm}\mathbb{E}\left[\sup_{0\leq t\leq T}e^{\frac{p}{2}\beta A_t}|Y^n_{t}|^{p}
		+\int_{0}^{T}e^{\frac{p}{2}\beta A_t}|Y^n_{t}|^{p}dA_t+\left(\int_{0}^{T}e^{\beta A_t}|Z^n_{t}|^{2}dt\right)^\frac{p}{2}+\left(\int_{0}^{T}e^{\beta A_t}d[M]_t\right)^\frac{p}{2}\right.\\
		&&\hspace{6cm}\left.+\left(\int_0^T\int_{\mathcal{U}}e^{\beta A_t}|V^n_t(e)|^2\mu(dt,de)\right)^\frac{p}{2}+|K^n_{T}|^p\right]\nonumber\\
		&\leq&\mathsf{K}_{p,\mathfrak{C},\beta}\mathbb{E}\left[e^{\frac{p}{2}\beta A_T}|\xi|^p+\left(\int_{0}^{T}e^{\beta A_t}\left|\frac{\mathfrak{f}(t)}{a_t}\right|^2dt\right)^\frac{p}{2}+\sup_{0\leq t\leq T}\left|e^{\frac{q}{2}\beta A_t}L_t^+\right|^{p}\right].
	\end{eqnarray*}
	Indeed, by applying Lemma \ref{Ito}, we get
	\begin{equation}\label{e01-n}
		\begin{split}
			&e^{\frac{p}{2}\beta A_{t \wedge \tau}}|Y^n_t|^p+\frac{p}{2}\beta \int_{t \wedge \tau}^\tau e^{\frac{p}{2}\beta A_s}|Y^n_s|^pdA_s+c(p)\int_{t \wedge \tau}^\tau e^{\frac{p}{2}\beta A_s}|Y^n_s|^{p-2}|Z^n_s|^2\mathds{1}_{\{Y^n_s\neq0\}}ds\\
			&+c(p)\int_{t \wedge \tau}^\tau e^{\frac{p}{2}\beta A_s}| Y^n_s|^{p-2}\mathds{1}_{\{Y^n_s\neq0\}}d[M^n]^c_s\\
			&\leq e^{\frac{p}{2}\beta A_T}|Y_\tau|^p+p\int_{t \wedge \tau}^\tau e^{\frac{p}{2}\beta A_s}|Y^n_s|^{p-1}\hat{Y}^n_s \mathfrak{f}(s)ds+p\int_{t \wedge \tau}^\tau e^{\frac{p}{2}\beta A_s}|Y^n_{s-}|^{p-1}\hat{Y}^n_{s-}dK^n_s\\
			&-p\int_{t \wedge \tau}^\tau e^{\frac{p}{2}\beta A_s}|Y^n_s|^{p-1}\hat{Y}^n_sZ^n_sdB_s-p\int_{t \wedge \tau}^\tau\int_{\mathcal{U}}e^{\frac{p}{2}\beta A_s}|Y^n_{s-}|^{p-1}\hat{Y}^n_{s-}V^n_s(e)\tilde{\mu}(ds,de)\\
			&-p\int_{t \wedge \tau}^\tau e^{\frac{p}{2}\beta A_s}|Y^n_{s-}|^{p-1}\hat{Y}^n_sZ^n_sdM^n_s\\
			&-\int_{t \wedge \tau}^\tau \int_{\mathcal{U}}e^{\frac{p}{2}\beta A_s}\left[|Y^n_{s-}+ V^n_s(e)|^{p}-|Y^n_{s-}|^p-p| Y^n_{s-}|^{p-1}\check{Y}^n_{s-} V^n_s(e)\right]\mu(ds,de)\\
			&-\sum_{t \wedge \tau < s \leq \tau } e^{\frac{p}{2}\beta A_s}\left[|\bar Y_{s-}+\Delta M^n _s|^{p}-|Y^n_{s-}|^p-p| Y^n_{s-}|^{p-1}\check{Y}^n_{s-} \Delta M^n_s\right].
		\end{split}
	\end{equation}
	By Hölder's
	inequality, we obtain
	\begin{equation}\label{eq1}
		\begin{split}
			&p\int_{t \wedge \tau}^\tau e^{\frac{p}{2}\beta A_s}| Y^n_{s}|^{p-1}\check{ Y}^n_{s}\mathfrak{f}(s)ds\\
			&\leq p\int_t^Te^{\frac{p}{2}\beta A_s}| Y^n_s|^{p-1}|\mathfrak{f}(s)| ds\\
			&=p\int_t^T\left(e^{\frac{p-1}{2}\beta A_s}|a_s|^{\frac{2(p-1)}{p}}|\bar Y_s|^{p-1}\right)\left(e^{\frac{\beta}{2}A_s}|a_s|^{\frac{2-p}{p}}\left|\frac{\mathfrak{f}(s)}{a_s}\right|\right)ds\\
			&\leq(p-1)\int_t^Te^{\frac{p}{2}\beta A_s}| Y^n_s|^{p}dA_s+\int_t^T\left(e^{\frac{p}{2}\beta A_s}\left|\frac{\mathfrak{f}(s)}{a_s}\right|^p\right)\left(|a_s|^{2-p}\right) ds\\
			&\leq(p-1)\int_t^Te^{\frac{p}{2}\beta A_s}| Y^n_s|^{p}dA_s+\left(\int_t^Te^{\beta A_s}\left|\frac{\mathfrak{f}(s)}{a_s}\right|^2ds\right)^\frac{p}{2}A_T^{\frac{2-p}{2}}.
		\end{split}
	\end{equation}
	Therefore, 
	\begin{equation}\label{finite1}
		\begin{split}
			\mathbb{E}\left[\int_0^T e^{\frac{p}{2}\beta A_s}| Y^n_{s}|^{p-1}\check{ Y}^n_{s}\mathfrak{f}(s)ds \right] <+\infty.
		\end{split}
	\end{equation}
	Moreover, it hold true that
	\begin{equation*}
		\int_{t \wedge \tau}^\tau e^{\frac{p}{2}\beta A_s}|Y^n_{s-}|^{p-1}\hat{Y}^n_s(Y^n_s-L_s)^-ds\leq\int_{t \wedge \tau}^\tau e^{\frac{p}{2}\beta A_s}|L^+_s|^{p-1}(Y^n_s-L_s)^-ds.
	\end{equation*}
	Hence, for each $\varrho>0$
	\begin{equation}\label{eq9-n}
		\begin{split}
			\int_{t \wedge \tau}^\tau e^{\frac{p}{2}\beta A_s}|Y^n_{s-}|^{p-1}\hat{Y}^n_{s-}dK^n_s
			&\leq\left(\int_0^T\left(e^{\frac{p}{2}\beta A_s}|L_s^+|^{p-1}\right)^\frac{p}{p-1}dK^n_s\right)^\frac{p-1}{p}.\left|K^n_T\right|^\frac{1}{p}\\
			&\leq\left(\int_0^T\left|e^{\frac{q}{2}\beta A_s}L_s^+\right|^{p}dK^n_s\right)^\frac{p-1}{p}.\left|K^n_T\right|^\frac{1}{p}\\
			&\leq\left(\sup_{0\leq t\leq T}\left|e^{\frac{q}{2}\beta A_t}L_t^+\right|^{p}\right)^\frac{p-1}{p}.\left|K^n_T\right|\\
			&\leq\frac{p-1}{p}\varrho^\frac{1}{p-1}\sup_{0\leq t\leq T}\left|e^{\frac{q}{2}\beta A_t}L_t^+\right|^{p}+\frac{1}{p\varrho}\left|K^n_T\right|^p.
		\end{split}
	\end{equation}
	Note also that
	$$
	\mathbb{E}\left[\left|K^n_T\right|^p\right] \leq \mathsf{K}_{n,p,T}\left(\mathbb{E}\left[\sup_{t\in[0,T]}  |{Y}^n_t|^p \right]+\mathbb{E}\left[\sup_{t\in[0,T]}  |L^+_t|^p \right]\right)<+\infty.
	$$
	Therefore, using assumption $(\mathcal{H}3)$-(ii), we derive
	\begin{equation}\label{Garcoffe}
		\mathbb{E}\left[\int_{0}^T e^{\frac{p}{2}\beta A_s}|Y^n_{s-}|^{p-1}\hat{Y}^n_{s-}dK^n_s\right]<+\infty.
	\end{equation}
	By the convexity of the function $x \mapsto |x|^p$ for $p \in (1,2)$, we obtain
	\begin{equation}\label{cv1}
		\begin{split}
			0 &\leq \int_{t \wedge \tau}^\tau \int_{\mathcal{U}}e^{\frac{p}{2}\beta A_s}\left[| Y^n_{s-}+ V^n_s(e)|^{p}-| Y^n_{s-}|^p-p| Y^n_{s-}|^{p-1}\check{ Y}^n_{s-} V^n_s(e)\right]\mu(ds,de)\\
			& \leq e^{\frac{p}{2}\beta A_\tau}|{Y}^n_\tau|^p+p\int_{t \wedge \tau}^\tau e^{\frac{p}{2}\beta A_s}| Y^n_s|^{p-1}\check{ Y}^n_s\mathfrak{f}(s)ds-p\int_t^Te^{\frac{p}{2}\beta A_s}|\bar Y^n_s|^{p-1}\check{ Y}^n_s  Z^n_sdB_s\\
			&+p\int_{t \wedge \tau}^\tau e^{\frac{p}{2}\beta A_s}|Y^n_{s-}|^{p-1}\hat{Y}^n_{s-}dK^n_s-p\int_{t \wedge \tau}^\tau \int_{\mathcal{U}}e^{\frac{p}{2}\beta A_s}| Y^n_{s-}|^{p-1}\check{ Y}^n_{s-} V^n_s(e)\tilde{\mu}(ds,de)\\
			&-\int_{t \wedge \tau}^\tau e^{\frac{p}{2}\beta A_s}| Y^n_{s-}|^{p-1}\check{ Y}^n_{s-}d M^n_s.
		\end{split}
	\end{equation}
	Consider a fundamental sequence of stopping times $\{\tau_m\}_{m \geq 1}$ associated with the following local martingale:
	\begin{equation*}
		\begin{split}
			&\int_0^\cdot e^{\frac{p}{2}\beta A_s}| Y^n_s|^{p-1}\check{Y}^n_s \bar Z^n_sdB_s+\int_{0}^\cdot \int_{\mathcal{U}}e^{\frac{p}{2}\beta A_s}|Y^n_{s-}|^{p-1}\check{\bar Y}_{s-}\bar V^n_s(e)\tilde{\mu}(ds,de)\\
			&+\int_0^\cdot e^{\frac{p}{2}\beta A_s}| Y^n_{s-}|^{p-1}\check{ Y}^n_{s-}dM^n_s.
		\end{split}
	\end{equation*}
	By taking expectations in \eqref{e01-n} with $\tau = \tau_m$, the local martingale term vanishes. Letting $m \to +\infty$ and applying the monotone and dominated convergence theorems, together with \eqref{finite1}, \eqref{Garcoffe}, and assumption $(\mathcal{H}1)$, we arrive at
	$$
	0 \leq \mathbb{E}\left[  \int_{0}^T \int_{\mathcal{U}}e^{\frac{p}{2}\beta A_s}\left[| Y^n_{s-}+ V_s(e)|^{p}-| Y^n_{s-}|^p-p| Y^n_{s-}|^{p-1}\check{ Y}^n_{s-} V^n_s(e)\right]\mu(ds,de)\right] <+\infty.
	$$
	Moreover, applying \cite[Lemma 3.67]{Jacod1979}, we also obtain the following
	$$
	0 \leq \mathbb{E}\left[  \int_{0}^T \int_{\mathcal{U}}e^{\frac{p}{2}\beta A_s}\left[| Y^n_{s-}+ V^n_s(e)|^{p}-|\bar Y_{s-}|^p-p|Y^n_{s-}|^{p-1}\check{Y}^n_{s-} V^n_s(e)\right] \lambda(de)ds \right] <+\infty.
	$$
	Applying Lemma 9 from \cite{kruse2016bsdes}, we have
	$$
	\begin{aligned}
		& \int_{t \wedge \tau}^\tau \int_{\mathcal{U}} e^{\frac{p}{2}\beta A_s} \left[\left|Y^n_{s-}+ V^n_s(e)\right|^p-\left| Y^n_{s-}\right|^p-p\left| Y^n_{s-}\right|^{p-1} \check{Y}^n_{s-}  V^n_s(e)\right] \mu(ds,de) \\
		& \quad \geq c(p) \int_{t \wedge \tau}^\tau \int_{\mathcal{U}} e^{\frac{p}{2}\beta A_s} \left| V^n_s(e)\right|^2\left(\left|Y^n_{s-}\right|^2 \vee\left| Y^n_{s-}+V^n_s(e)\right|^2\right)^{\frac{p-2}{2}}  \times \mathbf{1}_{\left| Y^n_{s-}\right| \vee\left|Y^n_{s-}+  V^n_s(e)\right| \neq 0} \mu(ds,de)
	\end{aligned}
	$$
	and
	$$
	\begin{aligned}
		& \sum_{t \wedge \tau<s \leq \tau} e^{\frac{p}{2}\beta A_s} \left[\left| Y^n_{s-}+\Delta M^n_s\right|^p-\left|\bar Y_{s-}\right|^p-p\left|\bar Y_{s-}\right|^{p-1} \check{Y}^n_{s-} \Delta M^n_s\right] \\
		& \quad \geq c(p) \sum_{t \wedge \tau<s \leq \tau} e^{\frac{p}{2}\beta A_s} \left|\Delta M^n_s\right|^2\left(\left| Y^n_{s-}\right|^2 \vee\left| Y^n_{s-}+\Delta M^n_s\right|^2\right)^{\frac{p-2}{2}} \mathbf{1}_{\left|\bar Y^n_{s-}\right| \vee\left| Y^n_{s-}+\Delta M^n_s\right| \neq 0}.
	\end{aligned}
	$$
	Coming back to \eqref{e01-n}, using using \eqref{eq1}, and the preceding estimates, we obtain
	\begin{equation}\label{start2}
		\begin{split}
			&e^{\frac{p}{2}\beta A_{t \wedge \tau}}|{Y}^n_{t \wedge \tau}|^p+\frac{p}{2}\beta\int_{t \wedge \tau}^ \tau e^{\frac{p}{2}\beta A_s}| Y^n_s|^pdA_s+{c(p)}\int_{t \wedge \tau}^\tau e^{\frac{p}{2}\beta A_s}|\bar Y_s|^{p-2}| Z^n_s|^2\mathds{1}_{\{Y_s\neq0\}}ds\\
			&+c(p)\int_{t \wedge \tau}^\tau e^{\frac{p}{2}\beta A_s}|Y^n_s|^{p-2}\mathds{1}_{\{Y^n_s\neq0\}}d[\bar M^n]^c_s\\
			&+{c(p)} \int_{t \wedge \tau}^\tau \int_{\mathcal{U}} e^{\frac{p}{2}\beta A_s} \left| V^n_s(e)\right|^2\left(\left| Y^n_{s-}\right|^2 \vee\left| Y^n_{s-}+ V^n_s(e)\right|^2\right)^{\frac{p-2}{2}}  \times \mathbf{1}_{\left|Y^n_{s-}\right| \vee\left| Y^n_{s-}+  V^n_s(e)\right| \neq 0} \mu(ds,de)\\
			&+c(p) \sum_{t \wedge \tau<s \leq \tau} e^{\frac{p}{2}\beta A_s} \left|\Delta  M^n_s\right|^2\left(\left|Y^n_{s-}\right|^2 \vee\left| Y^n_{s-}+\Delta M^n_s\right|^2\right)^{\frac{p-2}{2}} \mathbf{1}_{\left|Y^n_{s-}\right| \vee\left| Y^n_{s-}+\Delta M^n_s\right| \neq 0}\\
			&\leq e^{\frac{p}{2}\beta A_\tau}|{Y}^n_\tau|^p+\left(p-1\right)\int_{t \wedge \tau}^\tau e^{\frac{p}{2}\beta A_s}| Y^n_s|^{p}dA_s+\mathfrak{C}^{\frac{2-p}{2}}\left(\int_{t \wedge \tau}^\tau e^{\beta A_s}\left|\frac{\mathfrak{f}(s)}{a_s}\right|^2ds\right)^\frac{p}{2}\\
			&+{(p-1)}\varrho^\frac{1}{p-1}\sup_{0\leq t\leq T}\left|e^{\frac{q}{2}\beta A_t}L_t^+\right|^{p}+\frac{1}{\varrho}\left|K^n_T\right|^p-p\int_{t \wedge \tau}^\tau e^{\frac{p}{2}\beta A_s}| Y^n_s|^{p-1}\check{Y}^n_s Z^n_sdB_s\\
			& -p\int_{t \wedge \tau}^\tau \int_{\mathcal{U}}e^{\frac{p}{2}\beta A_s}|Y_{s-}|^{p-1}\check{Y}^n_{s-} V^n_s(e)\tilde{\mu}(ds,de)-p\int_{t \wedge \tau}^\tau e^{\frac{p}{2}\beta A_s}|Y^n_{s-}|^{p-1}\check{ Y}^n_{s-}d M^n_s.
		\end{split}
	\end{equation}
	Let $\{{\tau}_m\}_{m \geq 1}$ be a fundamental sequence for the local martingale term
	\begin{equation*}
		\begin{split}
			&\int_0^\cdot e^{\frac{p}{2}\beta A_s}| Y^n_s|^{p-1}\check{ Y}^n_s  Z^n_sdB_s+\int_{0}^\cdot \int_{\mathcal{U}}e^{\frac{p}{2}\beta A_s}|Y^n_{s-}|^{p-1}\check{Y}^n_{s-} V_s(e)\tilde{\mu}(ds,de)\\
			&+\int_{0}^\cdot \int_{\mathcal{U}}e^{\frac{p}{2}\beta A_s}\left[| Y^n_{s-}+ V^n_s(e)|^{p}-| Y^n_{s-}|^p-p| Y^n_{s-}|^{p-1}\check{ Y}^n_{s-}V^n_s(e)\right]\tilde{\mu}(ds,de).
		\end{split}
	\end{equation*}
	Next, taking $\tau = {\tau}_m$ then the expectation on both sides of the obtained estimation \eqref{start2} lead to the following estimates:
	\begin{equation}\label{start31}
		\begin{split}
			&\mathbb{E}\left[  e^{\frac{p}{2}\beta A_{t \wedge \tau}}|\bar \bar{Y}_{t \wedge \tau}|^p\right] +\left(\frac{p}{2}\beta-(p-1)\right)\mathbb{E}\left[\int_{t \wedge \tau}^ \tau e^{\frac{p}{2}\beta A_s}|\bar Y_s|^pdA_s\right]  \\
			+&{c(p)}\mathbb{E}\left[\int_{t \wedge \tau}^\tau e^{\frac{p}{2}\beta A_s}|\bar Y_s|^{p-2}|\bar Z_s|^2\mathds{1}_{\{Y_s\neq0\}}ds\right]+c(p)\mathbb{E}\left[\int_{t \wedge \tau}^\tau e^{\frac{p}{2}\beta A_s}|\bar Y_s|^{p-2}\mathds{1}_{\{\bar Y_s\neq0\}}d[M]^c_s\right] \\
			+&{c(p)} \mathbb{E}\left[\int_{t \wedge \tau}^\tau \int_{\mathcal{U}} e^{\frac{p}{2}\beta A_s} \left|\bar V_s(e)\right|^2\left(\left|\bar Y_{s-}\right|^2 \vee\left|\bar Y_{s-}+\psi_s(u)\right|^2\right)^{\frac{p-2}{2}}  \times \mathbf{1}_{\left|\bar Y_{s-}\right| \vee\left|\bar Y_{s-}+ \bar V_s(e)\right| \neq 0} \mu(ds,de)\right] \\
			+&c(p)\mathbb{E}\left[ \sum_{t \wedge \tau<s \leq \tau} e^{\frac{p}{2}\beta A_s} \left|\Delta M_s\right|^2\left(\left|\bar Y_{s-}\right|^2 \vee\left|\bar Y_{s-}+\Delta M_s\right|^2\right)^{\frac{p-2}{2}} \mathbf{1}_{\left|\bar Y_{s-}\right| \vee\left|\bar Y_{s-}+\Delta M_s\right| \neq 0}\right] \\
			&\leq \mathbb{E}\left[e^{\frac{p}{2}\beta A_\tau}|\bar{Y}_\tau|^p\right] +\mathfrak{C}^{\frac{2-p}{2}}\mathbb{E}\left[\left(\int_{t \wedge \tau}^\tau e^{\beta A_s}\left|\frac{\mathfrak{f}(s)}{a_s}\right|^2ds\right)^\frac{p}{2}\right]+(p-1)\varrho^\frac{1}{p-1}\mathbb{E}\left[ \sup_{0\leq t\leq T}\left|e^{\frac{q}{2}\beta A_t}L_t^+\right|^{p}\right]\\ 
			&~ +\frac{1}{\varrho}\mathbb{E}\left[\left|K^n_T\right|^p\right] .
		\end{split}
	\end{equation}
	By choosing $\beta > \frac{2}{p}(p-1)$ and using \eqref{finite1}, assumption $(\mathcal{H}1)$, together with the monotone and dominated convergence theorems, we can let $m \rightarrow +\infty$ in \eqref{start31} (for the sequence $\{\tau_m\}_{m \geq 1}$) to obtain
	\begin{equation}\label{start4}
		\begin{split}
			&\mathbb{E}\left[\int_{0}^T e^{\frac{p}{2}\beta A_s}| Y^n_s|^pdA_s\right] +\mathbb{E}\left[\int_{0}^T e^{\frac{p}{2}\beta A_s}|Y^n_s|^{p-2}|\bar Z_s|^2\mathds{1}_{\{Y^n_s\neq0\}}ds\right] \\
			&+\mathbb{E}\left[\int_{0}^T e^{\frac{p}{2}\beta A_s}| Y^n_s|^{p-2}\mathds{1}_{\{Y^n_s\neq0\}}d[ M^n]^c_s\right] \\
			&+ \mathbb{E}\left[\int_{0}^T \int_{\mathcal{U}} e^{\frac{p}{2}\beta A_s} \left| V^n_s(e)\right|^2\left(\left|Y^n_{s-}\right|^2 \vee\left| Y^n_{s-}+ V^n_s(e)\right|^2\right)^{\frac{p-2}{2}}  \times \mathbf{1}_{\left|Y^n_{s-}\right| \vee\left|Y^n_{s-}+  V^n_s(e)\right| \neq 0} \mu(ds,de)\right] \\
			&+\mathbb{E}\left[ \sum_{0<s \leq T} e^{\frac{p}{2}\beta A_s} \left| M^n_s\right|^2\left(\left| Y^n_{s-}\right|^2 \vee\left|Y^n_{s-}+\Delta M^n_s\right|^2\right)^{\frac{p-2}{2}} \mathbf{1}_{\left| Y^n_{s-}\right| \vee\left| Y^n_{s-}+\Delta M^n_s\right| \neq 0}\right] \\
			&\leq \mathsf{K}_{p,\mathfrak{C},\beta} \left(\mathbb{E}\left[e^{\frac{p}{2}\beta A_T}|{\xi}|^p\right] +\mathbb{E}\left[\left(\int_{0}^T e^{\beta A_s}\left|\frac{\mathfrak{f}(s)}{a_s}\right|^2ds\right)^\frac{p}{2}\right]+\mathbb{E}\left[ \sup_{0\leq t\leq T}\left|e^{\frac{q}{2}\beta A_t}L_t^+\right|^{p}\right] +\frac{1}{\varrho}\mathbb{E}\left[\left|K^n_T\right|^p\right] \right).
		\end{split}
	\end{equation}
	Let us now turn to \eqref{e01-n}. Using \eqref{eq1} and \eqref{eq9-n} with the same choice of $\beta$, we obtain
	\begin{equation*}\label{cv4}
		\begin{split}
			0 &\leq \mathbb{E}\left[ \int_{0}^T \int_{\mathcal{U}}e^{\frac{p}{2}\beta A_s}\left[| Y^n_{s-}+\bar V_s(e)|^{p}-|Y^n_{s-}|^p-p| Y^n_{s-}|^{p-1}\check{Y}^n_{s-} V^n_s(e)\right]\mu(ds,de)\right] \\
			& \leq \mathbb{E}\left[  e^{\frac{p}{2}\beta A_T}|{\xi}|^p\right]  +\mathfrak{C}^{\frac{2-p}{2}}\mathbb{E}\left[\left(\int_t^Te^{\beta A_s}\left|\frac{\mathfrak{f}(s)}{a_s}\right|^2ds\right)^\frac{p}{2}\right]+	(p-1)\varrho^\frac{1}{p-1}\mathbb{E}\left[ \sup_{0\leq t\leq T}\left|e^{\frac{q}{2}\beta A_t}L_t^+\right|^{p}\right] \\
			&  +\frac{1}{\varrho}\mathbb{E}\left[ \left|K^n_T\right|^p\right] +p\mathbb{E}\left[\sup_{t\in[0,T]}\left\{\left| \Gamma^n_t \right| +\left| \Theta^n_t\right| +\left| \Xi^n_t\right| \right\}\right], 
		\end{split}
	\end{equation*}
	where
	$$
	\begin{aligned}
		\Gamma^n_t & =\int_0^t e^{\frac{p}{2}\beta A_s}\left| Y^n_s\right|^{p-1} \check{Y}^n_s  Z^n_s d B_s \\
		\Theta^n_t & =\int_0^t e^{\frac{p}{2}\beta A_s}\left| Y^n_s\right|^{p-1} \check{Y}^n_s d M^n_s, \quad \Xi^n_t=\int_0^t e^{\frac{p}{2}\beta A_s}\left|Y^n_s\right|^{p-1} \check{ Y}^n_s \int_{\mathcal{U}} V^n_s(e) \tilde{\mu}(ds,de).
	\end{aligned}
	$$
	Since the integrands process $\ast$ in $\int_{0}^{\cdot} \ast \,\mu(ds,de)$ is predictable and, by definition of the solution, locally integrable with respect to $\mu$, we may take the predictable projection in the inequality above via conditional expectation using \cite[Theorem 3.15]{Jacod1979}. This yields
	\begin{equation}\label{cv5}
		\begin{split}
			0 &\leq \mathbb{E}\left[ \int_{0}^T \int_{\mathcal{U}}e^{\frac{p}{2}\beta A_s}\left[| Y^n_{s-}+\bar V_s(e)|^{p}-|Y^n_{s-}|^p-p| Y^n_{s-}|^{p-1}\check{Y}^n_{s-} V^n_s(e)\right]\lambda(de)ds\right] \\
			& \leq \mathbb{E}\left[  e^{\frac{p}{2}\beta A_T}|{\xi}|^p\right]  +\mathfrak{C}^{\frac{2-p}{2}}\mathbb{E}\left[\left(\int_t^Te^{\beta A_s}\left|\frac{\mathfrak{f}(s)}{a_s}\right|^2ds\right)^\frac{p}{2}\right]+	(p-1)\varrho^\frac{1}{p-1}\mathbb{E}\left[ \sup_{0\leq t\leq T}\left|e^{\frac{q}{2}\beta A_t}L_t^+\right|^{p}\right] \\
			& +\frac{1}{\varrho}\mathbb{E}\left[ \left|K^n_T\right|^p\right]+p\mathbb{E}\left[\sup_{t\in[0,T]}\left\{\left| \Gamma^n_t \right| +\left| \Theta^n_t\right| +\left| \Xi^n_t\right| \right\}\right].
		\end{split}
	\end{equation}
	Applying a convexity argument once more together with \eqref{cv5}, we take the supremum in \eqref{start2} followed by the expectation, and use the Burkholder-Davis-Gundy inequality (Theorem 48 in \cite[page 193]{Protter2004}) to obtain
	\begin{equation}\label{walo}
		\begin{split}
			\mathbb{E}\left[\sup_{t\in[0,T]} e^{\frac{p}{2}\beta A_t} |{Y}^n_t|^p \right]
			&\leq \mathbb{E}\left[  e^{\frac{p}{2}\beta A_T}|{\xi}|^p\right]  +\mathfrak{C}^{\frac{2-p}{2}}\mathbb{E}\left[\left(\int_t^Te^{\beta A_s}\left|\frac{\mathfrak{f}(s)}{a_s}\right|^2ds\right)^\frac{p}{2}\right] \\
			&+	(p-1)\varrho^\frac{1}{p-1}\mathbb{E}\left[ \sup_{0\leq t\leq T}\left|e^{\frac{q}{2}\beta A_t}L_t^+\right|^{p}\right]   +\frac{1}{\varrho}\mathbb{E}\left[ \left|K^n_T\right|^p\right]\\
			&+\mathsf{K}_p\mathbb{E}\left[\left(\left[ \Gamma^n\right]^{\frac{1}{2}}_T  +\left[  \Theta^n_t\right]^{\frac{1}{2}}_T  + \left[  \Xi^n \right]_T^{\frac{1}{2}}\right) \right].
		\end{split}
	\end{equation}
	\begin{itemize}
		\item The term $\left[ \Gamma^n\right]^{\frac{1}{2}}_T$ can be controlled as in \cite{Briand2003}:
		\begin{equation*}
			\begin{split}
				\mathsf{K}_p \mathbb{E}\left(  \left[ \Gamma^n \right]_T^{\frac{1}{2}}\right)   
				&\leq c \mathsf{K}_p\mathbb{E}\left[ \left(\int_0^Te^{p\beta A_s}| Y^n_s|^{2(p-1)}\mathds{1}_{\{ Y^n_s\neq 0\}}| Z^n_s|^2ds\right)^{\frac{1}{2}}\right] \\
				&\leq c \mathsf{K}_p \mathbb{E}\left[ \left(\sup_{0\leq t\leq T}e^{\frac{p}{2}\beta A_t}| Y^n_t|^{p}\int_0^Te^{\frac{p}{2}\beta A_s}| Y^n_s|^{p-2}\mathds{1}_{\{ Y^n_s\neq 0\}}| Z^n_s|^2ds\right)^{\frac{1}{2}}\right] \\
				&\leq \frac{1}{6}\mathbb{E}\left[\sup_{0\leq t\leq T}e^{\frac{p}{2}\beta A_t}| Y^n_t|^{p}\right]+\frac{3}{2} (c \mathsf{K}_p)^2\mathbb{E}\left[ \int_0^Te^{\frac{p}{2}\beta A_s}| Y^n_s|^{p-2}\mathds{1}_{\{ Y^n_s\neq 0\}}| Z^n_s|^2ds\right]. 
			\end{split}
		\end{equation*}

		\item The term $\left[ \Xi^n \right]^{\frac{1}{2}}_T$ can be controlled as in \cite{kruse2017lp}:
		\begin{equation*}
			\begin{split}
				\mathsf{K}_p \mathbb{E}\left(  \left[ \Xi^n \right]_T^{\frac{1}{2}}\right)  
				&\leq c\mathsf{K}_p\mathbb{E}\left[ \left(\int_0^T\int_{\mathcal{U}}e^{p\beta A_s}\left(| Y^n_{s-}|^{2}\vee| Y^n_{s-}+  V^n_s(e)|^2\right)^{p-1} \right. \right.\\ 
				&\left. \left.\qquad \qquad\qquad \times \mathds{1}_{\{| Y^n_{s-}|\vee| Y^n_{s-} +  V^n_s(e)|\neq 0\}}| V^n_s(e)|^2\mu(ds,de)\right)^{\frac{1}{2}}\right] \\
				&\leq\frac{1}{6}\mathbb{E}\left[\sup_{0\leq t\leq T}e^{\frac{p}{2}\beta A_t}| Y^n_t|^{p}\right]\\
				&+\frac{3}{2} (c\mathsf{K}_p)^2\mathbb{E}\left[ \int_0^T\int_{\mathcal{U}}e^{\frac{p}{2}\beta A_s}| V^n_s(e)|^2\left(| Y^n_{s-}|^{2}\vee| Y^n_{s-}+ V^n_s(e)|^2\right)^{\frac{p-2}{2}}\right.\\ 
				&\left.\qquad \qquad\qquad \times\mathds{1}_{\{| Y^n_{s-}|\vee| Y^n_{s-}+ V^n_s(e)|\neq 0\}}\mu(ds,de)\right] .
			\end{split}
		\end{equation*}

		\item The term $\left[ \Theta^n\right]^{\frac{1}{2}}_T$ can be controlled as in \cite{kruse2016bsdes}:
		\begin{equation*}
			\begin{split}
				&\mathsf{K}_p \mathbb{E}\left(  \left[  \Theta^n \right]_T^{\frac{1}{2}}\right) \\  
				&\leq\frac{1}{6}\mathbb{E}\left[\sup_{0\leq t\leq T}e^{\frac{p}{2}\beta A_t}| Y^n_t|^{p}\right]\\
				&+\frac{3}{2} (c \mathsf{K}_p)^2\mathbb{E}\left[\sum_{t \wedge \tau<s \leq \tau} e^{\frac{p}{2}\beta A_s} \left|\Delta  M^n_s\right|^2\left(\left| Y^n_{s-}\right|^2 \vee\left| Y^n_{s-}+\Delta M^n_s\right|^2\right)^{\frac{p-2}{2}} \mathbf{1}_{\left| Y^n_{s-}\right| \vee\left| Y^n_{s-}+\Delta  M^n_s\right| \neq 0}\right] \\
				&+(c\mathsf{K}_p)^2\mathbb{E}\left[\int_0^T e^{\frac{p}{2}\beta A_s}| Y^n_s|^{p-2}\mathds{1}_{\{ Y^n_s\neq0\}}d[ M^n]^c_s\right]
			\end{split}
		\end{equation*}	
	\end{itemize}
	Referring to \eqref{walo} and using \eqref{start4} along with the above estimations for the martingale parts, we obtain
	\begin{equation}\label{t}
		\begin{split}
			\mathbb{E}\left[\sup_{t\in[0,T]} e^{\frac{p}{2}\beta A_t} |{Y}^n_t|^p \right] 
			\leq & \mathsf{K}_{p,\mathfrak{C},\beta} \left(\mathbb{E}\left[e^{\frac{p}{2}\beta A_T}|{\xi}|^p\right] +\mathbb{E}\left[\left(\int_{0}^T e^{\beta A_s}\left|\frac{\mathfrak{f}(s)}{a_s}\right|^2ds\right)^\frac{p}{2}\right] \right.\\
			&\left. +\mathbb{E}\left[ \sup_{0\leq t\leq T}\left|e^{\frac{q}{2}\beta A_t}L_t^+\right|^{p}\right]   +\frac{1}{\varrho}\mathbb{E}\left[ \left|K^n_T\right|^p\right]\right).
		\end{split}
	\end{equation}
	As a result, combining \eqref{start4} with \eqref{t} yields, for every $\varrho > 0$, the estimate
	\begin{equation}\label{vd1}
		\begin{split}
			&\mathbb{E}\left[\sup_{t\in[0,T]} e^{\frac{p}{2}\beta A_t} |{Y}^n_t|^p \right]+\mathbb{E}\left[\int_{0}^{T}e^{\frac{p}{2}\beta A_s}| Y^n_{s}|^2dA_s \right] +\mathbb{E}\left[\left(\int_{0}^{T}e^{\beta A_s}| Z^n_{s}|^2ds\right)^{\frac{p}{2}}\right]\\
			&+	\mathbb{E}\left[\left(\int_0^T\int_{\mathcal{U}}e^{\beta A_t}| V^n_s(e)|^2\mu(ds,de)\right)^\frac{p}{2}\right]+\mathbb{E}\left[\left(\int_{0}^{T} e^{\beta A_s} d[M^n]_s\right)^{\frac{p}{2}} \right]\\
			&\leq  \mathsf{K}_{p,\mathfrak{C},\beta} \left(\mathbb{E}\left[e^{\frac{p}{2}\beta A_T}|{\xi}|^p\right] +\mathbb{E}\left[\left(\int_{0}^T e^{\beta A_s}\left|\frac{\mathfrak{f}(s)}{a_s}\right|^2ds\right)^\frac{p}{2}\right] +\mathbb{E}\left[ \sup_{0\leq t\leq T}\left|e^{\frac{q}{2}\beta A_t}L_t^+\right|^{p}\right] +\frac{1}{\varrho}\mathbb{E}\left[\left|K^n_T\right|^p\right] 	 \right).
		\end{split}
	\end{equation}
	By using the basic inequality 
	\begin{equation}\label{basic}
		\left(\sum_{i=1}^n|X_i|\right)^p\leq n^p\sum_{i=1}^n|X_i|^p\qquad\forall (n,p)\in\mathbb{N}^\ast\times]0,+\infty[,
	\end{equation}
	we get
	\begin{equation}\label{e007-n}
		\begin{split}
			\mathbb{E}\left|K^n_{T}\right|^p
			&\leq5^{p}\mathbb{E}\left[\sup_{0\leq t\leq T}e^{\frac{p}{2}\beta A_t}|Y^n_{t}|^p+e^{\frac{p}{2}\beta A_T}|\xi|^p+\left|\int_{0}^{T}\mathfrak{f}(s)ds\right|^{p}+\left|\int_{0}^{T}Z^n_{s}dB_{s}\right|^{p}\right.\\
			&\left.\qquad\qquad+\left|\int_0^T\displaystyle\int_{\mathcal{U}}V^n_s(e)\tilde{\mu}(ds,de)\right|^{p}+\left|\int_{0}^{T}dM^n_{s}\right|^{p}\right].
		\end{split}
	\end{equation}
	By the Burkholder-Davis-Gundy inequality, there exists a universal non-negative constant $c$ such that
	\begin{equation*}
		\mathbb{E}\left[ \sup_{0\leq t\leq T}\left|\int_0^tZ^n_sdB_s\right|^p\right] \leq c\mathbb{E}\left[ \left(\int_0^Te^{\beta A_s}|Z^n_s|^2ds\right)^{\frac{p}{2}}\right] 
	\end{equation*}
	and
	\begin{equation*}
		\mathbb{E}\left[ \sup_{0\leq t\leq T}\left|\int_0^t\int_{\mathcal{U}}V^n_s(e)\tilde{\mu}(ds,de)\right|^p\right] \leq c\mathbb{E}\left[ \left(\int_0^T\int_{\mathcal{U}}e^{\beta A_s}|V^n_s(e)|^2\mu(ds,de)\right)^{\frac{p}{2}}\right] 
	\end{equation*}
	and
	\begin{equation*}
		\mathbb{E}\left[ \sup_{0\leq t\leq T}\left|\int_0^t dM^n_s\right|^p\right] \leq c\mathbb{E}\left[ \left(\int_0^Te^{\beta A_s} d[M^n]_s\right)^{\frac{p}{2}}\right] .
	\end{equation*}
	By choosing a large $\varrho$ such that $\varrho>5^p \mathsf{K}_{p,\mathfrak{C},\beta}$ in \eqref{vd1}, we obtain
	\begin{equation*}\label{vd11}
		\begin{split}
			&\mathbb{E}\left[\sup_{t\in[0,T]} e^{\frac{p}{2}\beta A_t} |{Y}^n_t|^p \right] +\mathbb{E}\left[\int_{0}^{T}e^{\frac{p}{2}\beta A_s}| Y^n_{s}|^2dA_s \right] +\mathbb{E}\left[\left(\int_{0}^{T}e^{\beta A_s}| Z^n_{s}|^2ds\right)^{\frac{p}{2}}\right]\\
			&+	\mathbb{E}\left[\left(\int_0^T\int_{\mathcal{U}}e^{\beta A_t}| V^n_s(e)|^2\mu(ds,de)\right)^\frac{p}{2}\right]+\mathbb{E}\left[\left(\int_{0}^{T} e^{\beta A_s} d[M^n]_s\right)^{\frac{p}{2}} \right]\\
			&\leq  \mathsf{K}_{p,\mathfrak{C},\beta} \left(\mathbb{E}\left[e^{\frac{p}{2}\beta A_T}|{\xi}|^p\right] +\mathbb{E}\left[\left(\int_{0}^T e^{\beta A_s}\left|\frac{\mathfrak{f}(s)}{a_s}\right|^2ds\right)^\frac{p}{2}\right]	+\mathbb{E}\left[ \sup_{0\leq t\leq T}\left|e^{\frac{q}{2}\beta A_t}L_t^+\right|^{p}\right] \right)
		\end{split}
	\end{equation*}
	and then
	\begin{equation*}\label{e07-n}
		\mathbb{E}\left[ \left|K^n_{T}\right|^p\right] 
		\leq \mathsf{K}_{p,\mathfrak{C},\beta} \mathbb{E}\left[e^{\frac{p}{2}\beta A_T}|\xi|^p+\left(\int_0^Te^{\beta A_s}\left|\frac{\mathfrak{f}(s)}{a_s}\right|^2ds\right)^\frac{p}{2}+\sup_{0\leq t\leq T}\left|e^{\frac{q}{2}\beta A_t}L_t^+\right|^{p}\right].
	\end{equation*}
	
	~\\
	\textbf{Step 2:} There exists a progressively measurable process $Y$ such that 
	\begin{itemize}
		\item $Y \geq L$, and
		\item $\displaystyle \mathbb{E}\Bigl[\,\sup_{0 \leq t \leq T} e^{p\beta A_t} |(Y^n_t - L_t)^-|^p \Bigr] \xrightarrow[n \to +\infty]{} 0$.
	\end{itemize}
	
	~\\
	Indeed, from Proposition \ref{com}, we deduce that $Y^{n+1}_t \geq Y^n_t$ for each $n\in\mathbb{N}$. Hence, there exists a progressively measurable process $Y$ such that, for any $t \in [0,T]$, $Y_t:=\lim\limits_{n\rightarrow +\infty} Y^n_t$. Moreover, by Fatou's lemma, we have
	$
	\mathbb{E}\left[|Y_t|^p\right] <+\infty
	$
	since $\sup_{n \geq 0}\left\{\mathbb{E}\left[\sup_{t\in[0,T]} e^{\frac{p}{2}\beta A_t} |{Y}^n_t|^p \right] \right\} \leq  \mathsf{K}_{p,\mathfrak{C},\beta}$ from \textbf{Step 1}.\\
	On the other hand, from the previous step, we also know that $\sup_{n \geq 0}	\mathbb{E}\left[ \left|K^n_{T}\right|^p\right]\leq  \mathsf{K}_{p,\mathfrak{C},\beta}$. Therefore, taking the limit as $n \rightarrow +\infty$, we deduce that
	$$
	\mathbb{E}\left[\int_{0}^{T}\left(Y_s-L_s\right)^- ds\right]=0,
	$$
	and hence $\mathbb{P}$-a.s., $Y_t \geq L_t$ for any $t \in [0,T)$. As $\xi \leq L_T$ and $Y_T=\lim\limits_{n\rightarrow +\infty} Y^n_T=\xi$, it follows that $Y \geq L$. Next, if we denote by
	${}^p\!X$ the predictable projection of any process $X$ (see, e.g., \cite[Ch. V. Sec 1]{HeWangYan1992}), we have ${}^p\!Y^n \nearrow {}^p\!Y$ and  ${}^p\!Y \geq {}^p\!L$. \\
	For any $n \in \mathbb{N}$, the jump times of the process $\int_{0}^{\cdot} \int_{\mathcal{U}} V^n_s(e)\tilde{\mu}(ds,de)$ are totally
	inaccessible. It follows that the jump times of $Y^n$ are also totally inaccessible. Thus, for any predictable stopping time $\tau$, we have $\Delta Y^n_{\tau}=0$, and hence the predictable
	projection of $Y^n$ is the left-limited process $Y^n_-$, i.e., ${}^p\!Y^n=Y^n_{-}$. Similarly, by assumption, ${}^p\!L=L_{-}$. Therefore, we have proved that $Y^n_{-}={}^p\!Y^n \nearrow {}^p\!Y \geq {}^p\!L =L_{-}$. It follows that $(Y^n_{-}-L_-)^- \searrow ({}^p\! Y-L_{-})^- =0$.\\
	Consequently, by the generalized version of Dini's theorem (see, e.g., \cite[page 202]{DellacherieMeyer1980}), we deduce that $\sup_{0\leq t\leq T} e^{p\beta A_t}(Y^n_{-}-L_-)^- \searrow 0$ $\mathbb{P}$-a.s. as $n \rightarrow +\infty$. Furthermore, we have $\sup_{n \in \mathbb{N}}\sup_{0\leq t\leq T}e^{p\beta A_t}\left(Y^n_t-L_t\right)^- \leq  \sup_{0\leq t\leq T} e^{p\beta A_t}|Y^0_t|+ \sup_{0\leq t\leq T}e^{p\beta A_t}|L^+_t|$ a.s. Therefore, the dominated convergence theorem implies
	$$
	\lim\limits_{n \rightarrow +\infty}\mathbb{E}\left[\sup\limits_{0\leq t\leq T}e^{p\beta A_t}|(Y^n_t-L_t)^-|^p\right]=0.
	$$
	
	~\\
	\textbf{Step 3:} There exists an $\mathbb{F}$-adapted process $(Y_t,Z_t,V_t,M_t,K_t)_{t\leq T}$ such that
	\begin{eqnarray*}
		&&\|Y^{n}-Y\|_{\mathfrak{B}^p_\beta}^{p}+\|Z^{n}-Z\|_{\mathcal{H}^p_\beta}^{p}+\|V^n-V\|_{\mathfrak{L}^p_{\mu,\beta}}^p+\|M^{n}-M\|_{\mathcal{M}^p_\beta}^{p}
		+\|K^{n}-K\|_{\mathcal{S}^p}^{p}\xrightarrow[n\to +\infty]{}0.
	\end{eqnarray*}
	
	Indeed, let $\Re^{n,m}=\Re^n-\Re^m$ for each $n\geq m\geq0$ and for $\Re\in\{Y,Z,V,M,K\}$. Once again, Lemma \ref{Ito} implies
	\begin{equation}\label{eq6}
		\begin{split}
			&e^{\frac{p}{2}\beta A_t}|Y^{n,m}_t|^p+\frac{p}{2}\beta\int_t^Te^{\frac{p}{2}\beta A_s}|Y^{n,m}_s|^pdA_s+c(p)\int_t^Te^{\frac{p}{2}\beta A_s}|Y^{n,m}_s|^{p-2}|Z^{n,m}_s|^2\mathds{1}_{\{Y^n_s\neq Y^m_s\}}ds\nonumber\\
			&+c(p)\int_t^T\int_{\mathcal{U}}e^{\frac{p}{2}\beta A_s}|V^{n,m}_s(e)|^2\left(|Y^{n,m}_{s-}|^{2}\vee|Y^{n,m}_{s}|^2\right)^{\frac{p-2}{2}}\mathds{1}_{\{|Y^{n,m}_{s-}|\vee|Y^{n,m}_{s-}+V^{n,m}_s(e)|\neq 0\}}\mu(ds,de)\nonumber\\
			&+c(p)\int_{t \wedge \tau}^\tau e^{\frac{p}{2}\beta A_s}| Y^{n,m}_s|^{p-2}\mathds{1}_{\{ Y^{n,m}_{s}\neq0\}}d[M^{n,m}]^c_s \nonumber \\
			&\leq p\int_t^Te^{\frac{p}{2}\beta A_s}|Y^{n,m}_{s-}|^{p-1}\check{Y}^{n,m}_{s-} dK^{n,m}_s-p\int_t^Te^{\frac{p}{2}\beta A_s}|Y^{n,m}_{s}|^{p-1}\check{Y}^{n,m}_{s}Z^{n,m}_{s}dB_s \nonumber \\
			&-p\int_t^T\int_{\mathcal{U}}e^{\frac{p}{2}\beta A_s}|Y^{n,m}_{s-}|^{p-1}\check{Y}^{n,m}_{s} V^{n,m}_{s}(e)\tilde{\mu}(ds,de)
			-p\int_t^Te^{\frac{p}{2}\beta A_s}|Y^{n,m}_{s-}|^{p-1}\check{Y}^{n,m}_{s} dM^{m,n}_s.
		\end{split}
	\end{equation}
	By performing the same arguments as in \textbf{Step 1}, we obtain the existence of a constant $K_{p,\mathfrak{C},\beta}$ such that 
	\begin{equation}\label{vd111}
		\begin{split}
			&\mathbb{E}\left[\sup_{t\in[0,T]} e^{\frac{p}{2}\beta A_t} |{Y}^{n,m}_t|^p \right]+\mathbb{E}\left[\int_{0}^{T}e^{\frac{p}{2}\beta A_s}| Y^n_{s}|^2dA_s \right] +\mathbb{E}\left[\left(\int_{0}^{T}e^{\beta A_s}| Z^{n,m}_{s}|^2ds\right)^{\frac{p}{2}}\right]\\
			&+	\mathbb{E}\left[\left(\int_0^T\int_{\mathcal{U}}e^{\beta A_t}| V^{n,m}_s(e)|^2\mu(ds,de)\right)^\frac{p}{2}\right]+\mathbb{E}\left[\left(\int_{0}^{T} e^{\beta A_s} d[M^{n,m}]_s\right)^{\frac{p}{2}} \right]\\
			& \leq \mathsf{K}_{p,\mathfrak{C},\beta}  \mathbb{E}\left[ \left(\sup_{0\leq t\leq T}e^{p\beta A_t}|(Y^m_t-L_t)^-|^{p}\right)^\frac{p-1}{p}\right] \mathbb{E}\left[ \left(|K^n_T|^p\right)^\frac{1}{p}\right] \\
			&+\mathsf{K}_{p,\mathfrak{C},\beta} \mathbb{E}\left[ \left(\sup_{0\leq t\leq T}e^{p\beta A_t}|(Y^n_t-L_t)^-|^{p}\right)^\frac{p-1}{p}\right] \mathbb{E}\left[ \left(|K^m_T|^p\right)^\frac{1}{p}\right] \\
			& \leq \mathsf{K}_{p,\mathfrak{C},\beta}\left(\mathbb{E}\left(\sup_{0\leq t\leq T}e^{p\beta A_t}|(Y^m_t-L_t)^-|^{p}\right)^\frac{p-1}{p}+\mathbb{E}\left(\sup_{0\leq t\leq T}e^{p\beta A_t}|(Y^n_t-L_t)^-|^{p}\right)^\frac{p-1}{p}\right)\\
			&\xrightarrow[n,m\to +\infty]{}0.
		\end{split}
	\end{equation}
	Thus, $\left\lbrace Y^n\right\rbrace _{n\geq0}$ is a Cauchy sequence in $\mathcal{S}^{p,A}_\beta \cap \mathcal{S}^{p}_\beta$. Since $Y^m\nearrow Y$, we have $$\mathbb{E}\left[\sup_{t\in[0,T]} e^{\frac{p}{2}\beta A_t} |{Y}^{n}_t-Y_t|^p \right]\xrightarrow[n\to +\infty]{}0 ~\textnormal{ and  }~\mathbb{E}\int_0^Te^{\frac{p}{2}\beta A_s}|Y^n_s-Y_s|^pdA_s\xrightarrow[n\to +\infty]{}0,$$ 
	so $Y\in\mathcal{S}^{p,A}_\beta \cap \mathcal{S}^{p}_\beta$. \\
	We also obtain that $\left\lbrace (Z^n,V^n,M^n)\right\rbrace_{n\geq0}$ is a Cauchy sequence of processes in $\mathcal{H}^p_\beta\times\mathfrak{L}^p_{\mu,\beta} \times \mathcal{M}^p_\beta$. Hence, there exists a triplet of processes $(Z,V,M)$ such that the sequences $\left\lbrace Z^n\right\rbrace _{n\geq0}$, $\left\lbrace V^n\right\rbrace _{n\geq0}$, and $\left\lbrace M^n\right\rbrace _{n\geq0}$ converge to $Z\in\mathcal{H}^p_\beta$, $V\in\mathfrak{L}^p_{\mu,\beta}$, and $M\in\mathcal{M}^p_\beta$, respectively.\\
	To conclude, from (\ref{penalized}), we have
	$$
	K^n_t=Y^n_t-Y^n_0+\int_0^t \mathfrak{f}(s)ds-\int_0^tZ^n_sdB_s-\int_0^t\int_{\mathcal{U}}V^n_s(e)\tilde{\mu}(ds,de)-\int_0^t dM^n_s .
	$$
	Then
	$$
	\mathbb{E}\left[\sup_{0\leq t\leq T}|K_{t}^{n}-K_{t}^m|^{p}\right]\xrightarrow[n,m\to +\infty]{}0.
	$$
	It follows that $\left\lbrace K^n\right\rbrace_{n\geq0}$ is a Cauchy sequence in $\mathcal{S}^p$. Hence, there exists an $\mathbb{F}$-adapted, non-decreasing, continuous process $(K_t)_{t \leq T}$ with $K_0 = 0$ such that $\mathbb{E}\left[\sup_{0\leq t\leq T}|K_{t}^{n}-K_{t}|^{p}\right]\xrightarrow[n\to +\infty]{}0$.\\
	Finally, passing to the limit in $\mathbb{L}^p$ term by term in \eqref{penalized}, we obtain
	$$
	Y_{t}=\xi +\int_{t}^{T} \mathfrak{f}(s)ds+(K_T-K_t)-\int_{t}^{T}Z_{s}dB_{s}
	-\int_t^T\int_{\mathcal{U}}V_s(e)\tilde{\mu}(ds,de)-\int_{t}^{T} dM_s
	$$
	with $Y \geq L$, and $Y$ is RCLL since $K$ is continuous. It remains to show the Skorokhod condition $\int_{0}^{T}\left(Y_s-L_s\right)dK_s=0$.
	
	~\\
	\textbf{Step 4:} The limiting process $(K_t)_{t\leq T}$ verifies the Skorokhod condition $\int_{0}^{T}\left(Y_s-L_s\right)dK_s=0$. 
	
	This follows by applying the same reasoning as in Step~6 of the proof of \cite[Theorem~1.2(a)]{HMSOUK}. For this reason, we omit the detailed proof.
	
	\paragraph*{Part 2: General case of the driver $f$}\emph{}\\
	Let $(y,z,v,m) \in \mathcal{S}^{p,A}_\beta \times \mathcal{H}^p_\beta \times \mathfrak{L}^p_{\mu,\beta} \times \mathcal{M}^p_\beta$, and define $(Y,Z,V,M) = \Psi(y,z,v,m)$, where $(Y,Z,V,M,K)$ denotes the unique $\mathbb{L}^p$-solution of the RBSDE \eqref{basic RBSDE thm}.  
	For another element $(y',z',v',m') \in \mathcal{S}^{p,A}_\beta \times \mathcal{H}^p_\beta \times \mathfrak{L}^p_{\mu,\beta} \times \mathcal{M}^p_\beta$, we similarly set $(Y',Z',V',M') = \Psi(y',z',v',m')$, where $(Y',Z',V',M',K')$ is the unique $\mathbb{L}^p$-solution of \eqref{basic RBSDE thm} with parameters $(\xi, f(\cdot,y',z',v'), L)$.  \\
	We denote $\bar{\Re} = \Re - \Re'$ for $\Re \in \{Y,Z,V,K,Y',Z',V',K'\}$ and define $\delta f_t = f(t,y_t,z_t,v_t) - f(t,y'_t,z'_t,v'_t)$.  
	Our goal is to show that the mapping $\Psi$ is a strict contraction on $\mathcal{S}^{p,A}_\beta \times \mathcal{H}^p_\beta \times \mathfrak{L}^p_{\mu,\beta}  \times \mathcal{M}^p_\beta$, equipped with the norm
	$$
	\|(Y,Z,V,M)\|_{\mathcal{S}^{p,A}_\beta \times \mathcal{H}^p_\beta \times \mathfrak{L}^p_{\mu,\beta}}^p
	:= \|Y\|_{\mathcal{S}^{p,A}_\beta}^p + \|Z\|_{\mathcal{H}^p_\beta}^p + \|V\|_{\mathfrak{L}^p_{\mu,\beta}}^p + \|M\|_{\mathcal{M}^p_\beta}^p.
	$$
	By applying Lemma \ref{Ito} and using arguments that are now standard in this framework, we obtain:
	\begin{equation}\label{startF}
		\begin{split}
			&e^{\frac{p}{2}\beta A_{t \wedge \tau}}| \bar{Y}_{t \wedge \tau}|^p+\frac{p}{2}\beta\int_{t \wedge \tau}^ \tau e^{\frac{p}{2}\beta A_s}|\bar Y_s|^pdA_s+c(p)\int_{t \wedge \tau}^\tau e^{\frac{p}{2}\beta A_s}|\bar Y_s|^{p-2}|\bar Z_s|^2\mathds{1}_{\{\bar Y_s\neq0\}}ds\\
			&+c(p)\int_{t \wedge \tau}^\tau e^{\frac{p}{2}\beta A_s}|\bar Y_s|^{p-2}\mathds{1}_{\{\bar Y_s\neq0\}}d[\bar M]^c_s\\
			&\leq e^{\frac{p}{2}\beta A_\tau}|\bar{Y}_\tau|^p+p\int_t^Te^{\frac{p}{2}\beta A_s}|\bar{Y}_s|^{p-1}\check{\bar{Y}}_{s}\delta f_sds+p\int_t^Te^{\frac{p}{2}\beta A_s}|\bar{Y}_{s}|^{p-1}\check{\bar{Y}}_{s}d\bar{K}_s\\
			&-p\int_t^Te^{\frac{p}{2}\beta A_s}|\bar Y_s|^{p-1}\check{\bar Y}_s \bar Z_sdB_s-p\int_{t \wedge \tau}^\tau \int_{\mathcal{U}}e^{\frac{p}{2}\beta A_s}|\bar Y_{s-}|^{p-1}\check{\bar Y}_{s-}\bar V_s(e)\tilde{\mu}(ds,de)\\
			&-\int_{t \wedge \tau}^\tau \int_{\mathcal{U}}e^{\frac{p}{2}\beta A_s}\left[|\bar Y_{s-}+\bar V_s(e)|^{p}-|\bar Y_{s-}|^p-p|\bar Y_{s-}|^{p-1}\check{\bar Y}_{s-}\bar V_s(e)\right]\mu(ds,de)\\
			&-p\int_t^Te^{\frac{p}{2}\beta A_s}|\bar Y_{s-}|^{p-1}\check{\bar Y}_{s-}d\bar M_s-\sum_{t \wedge \tau < s \leq \tau } e^{\frac{p}{2}\beta A_s}\left[|\bar Y_{s-}+\Delta \bar M _s|^{p}-|\bar Y_{s-}|^p-p|\bar Y_{s-}|^{p-1}\check{Y}_{s-} \Delta \bar M_s\right].
		\end{split}
	\end{equation}
	Note that, for any square-integrable function $t \ni [0,T] \mapsto h(t)$, by Jensen's inequality, we have
	\begin{equation}\label{INSEA}
		\int_{0}^{t} e^{\frac{p}{2}\beta A_s} |h(s)|^{p} ds \leq t^{\frac{2-p}{2}} \left(\int_{0}^{t} e^{\beta A_s} |h(s)|^2 ds\right)^{\frac{p}{2}} ,\quad t \in [0,T],
	\end{equation}
	Moreover, by H\"older's inequality we have
	\begin{equation}\label{SKS1}
		\begin{split}
			p\int_{t \wedge \tau}^\tau e^{\frac{p}{2}\beta A_s}|\bar Y_s|^{p-1}| \bar y_{s}| \theta_s ds &\leq p \left(\int_{0}^T e^{\frac{p}{2}\beta A_s} |\bar{Y}_s|^p dA_s\right)^{\frac{p-1}{p}}\left(\int_{0}^T e^{\frac{p}{2}\beta A_s}| \bar y_{s}|^pds\right)^{\frac{1}{p}}
		\end{split}
	\end{equation}
	and
	\begin{equation}\label{SKS3}
		\begin{split}
			p\int_{t \wedge \tau}^\tau e^{\frac{p}{2}\beta A_s}|\bar Y_s|^{p-1}| \bar z_{s}| \eta_s ds &\leq p \left(\int_{0}^T e^{\frac{p}{2}\beta A_s} |\bar{Y}_s|^p dA_s\right)^{\frac{p-1}{p}}\left(\int_{0}^T e^{\frac{p}{2}\beta A_s}| \bar z_{s}|^pds\right)^{\frac{1}{p}}
		\end{split}
	\end{equation}
	and
	\begin{equation}\label{SKS2}
		\begin{split}
			p\int_{t \wedge \tau}^\tau e^{\frac{p}{2}\beta A_s}|\bar Y_s|^{p-1}\| \bar v_{s}\|_{\mathbb{L}^1_\lambda+\mathbb{L}^2_\lambda} \eta_s ds &\leq p \left(\int_{0}^T e^{\frac{p}{2}\beta A_s} |\bar{Y}_s|^p dA_s\right)^{\frac{p-1}{p}}\left(\int_{0}^T e^{\frac{p}{2}\beta A_s}\| \bar v_{s}\|^p_{\mathbb{L}^1_\lambda+\mathbb{L}^2_\lambda}ds\right)^{\frac{1}{p}}
		\end{split}
	\end{equation}
	So, according to the stochastic Lipschitz condition $(\mathcal{H}2)$-(ii)-(iii) on $f$, by using estimations \eqref{INSEA}, \eqref{SKS1}, \eqref{SKS2}, \eqref{SKS3} and Young's inequality, we have for any $\varrho>0$
	\begin{equation*}
		\begin{split}
			&p\int_t^Te^{\frac{p}{2}\beta A_s}|\bar Y_{s}|^{p-1}\hat{\bar Y}_{s}\delta f_sds\nonumber\\
			&\leq p\int_t^Te^{\frac{p}{2}\beta A_s}|\bar Y_{s}|^{p-1}(\theta_s|\bar y_s|+\gamma_s |\bar z_s|+\eta_s\|\bar v_{s}\|_{\mathbb{L}^1_\lambda+  \mathbb{L}^2_\lambda})ds\nonumber\\
			&\leq p\left(\int_t^Te^{\frac{p}{2}\beta A_s}\theta_s|\bar Y_{s}|^{p}ds\right)^\frac{p-1}{p}\left(\int_t^Te^{\frac{p}{2}\beta A_s}\theta_s|\bar y_{s}|^{p}ds\right)^\frac{1}{p}\\
			&+p\left(\int_t^Te^{\frac{p}{2}\beta A_s}|\gamma_s|^q|\bar Y_{s}|^{p}ds\right)^\frac{p-1}{p}\left(\int_t^Te^{\frac{p}{2}\beta A_s}|\bar z_{s}|^{p}ds\right)^\frac{1}{p}\\
			&+p\left(\int_t^Te^{\frac{p}{2}\beta A_s}|\eta_s|^q|\bar Y_{s}|^{p}ds\right)^\frac{p-1}{p}\left(\int_t^Te^{\frac{p}{2}\beta A_s}\|\bar v_{s}\|_{\mathbb{L}^1_\lambda+  \mathbb{L}^2_\lambda}^{p}ds\right)^\frac{1}{p}\nonumber\\
			&\leq 3(p-1)\varrho^\frac{p-1}{p^2}\int_t^Te^{\frac{p}{2}\beta A_s}|\bar Y_{s}|^{p}dA_s+\frac{1}{\varrho}\left(\int_t^Te^{\frac{p}{2}\beta A_s}|\bar y_{s}|^{p}dA_s+\left(\int_t^Te^{\beta A_s}|\bar z_{s}|^{2}ds\right)^\frac{p}{2}\right.\\
			&\qquad \qquad \qquad\left. +\int_t^Te^{\beta A_s}\|\bar v_{s}\|_{\mathbb{L}^1_\lambda+  \mathbb{L}^2_\lambda}^{p}ds\right).\nonumber\\
			&\leq 3(p-1)\varrho^\frac{p-1}{p^2}\int_t^Te^{\frac{p}{2}\beta A_s}|\bar Y_{s}|^{p}dA_s+\frac{1 \vee T^{\frac{p}{2}-1}}{\varrho}\left(\int_t^Te^{\frac{p}{2}\beta A_s}|\bar y_{s}|^{p}dA_s+\left(\int_t^Te^{\beta A_s}|\bar z_{s}|^{2}ds\right)^\frac{p}{2}\right.\\
			&\qquad \qquad \qquad\left.+\int_t^Te^{\beta A_s}\|\bar v_{s}\|_{\mathbb{L}^1_\lambda+  \mathbb{L}^2_\lambda}^{p}ds+\left(\int_t^Te^{\beta A_s}d [\bar{m}]_s\right)^{\frac{p}{2}}\right).
		\end{split}
	\end{equation*}
	The term $\int_t^T e^{\beta A_s} \|\bar{v}_{s}\|_{\mathbb{L}^1_\lambda + \mathbb{L}^2_\lambda}^{p} ds$ can be controlled by means of \cite[Lemma 2]{kruse2017lp}. More precisely, by applying the previous mentioned Lemma, we may derive the following estimation
	\begin{equation}\label{sK}
		\mathbb{E}\left[\int_{t}^{T} e^{\frac{p}{2}\beta A_s} \| \bar v_{s}\|_{\mathbb{L}^1_\lambda+\mathbb{L}^2_\lambda}^{p} ds \right] \leq K_{p,T} \mathbb{E}\left[\left(\int_{t}^{T} \int_{\mathcal{U}} e^{\beta A_s} |\bar v_s(e)|^2 \mu(ds,de) \right)^{\frac{p}{2}}\right].
	\end{equation}
	Moreover, we have $\int_t^T e^{\frac{p}{2} \beta A_s} |\bar{Y}_{s}|^{p-1} \hat{\bar{Y}}_{s-} \, d\bar{K}_s \leq 0$, since
	\begin{equation*}
		\begin{split}
			&\int_t^T e^{\frac{p}{2} \beta A_s} |\bar{Y}_{s}|^{p-1} \hat{\bar{Y}}_{s} \, d\bar{K}_s\\
			&\leq \int_t^T e^{\frac{p}{2} \beta A_s} |\bar{Y}_s|^{p-2} \mathds{1}_{\{\bar{Y}_s \neq 0\}} (Y_s - L_s) \, dK_s
			+ \int_t^T e^{\frac{p}{2} \beta A_s} |\bar{Y}_s|^{p-2} \mathds{1}_{\{\bar{Y}_s \neq 0\}} (Y'_s - L_s) \, dK'_s \\
			&= 0,
		\end{split}
	\end{equation*}
	where we used the facts that $dK_s = \mathds{1}_{\{Y_s = L_s\}} dK_s$ and $dK'_s = \mathds{1}_{\{Y'_s = L_s\}} dK'_s$.  \\
	Returning to (\ref{startF}) and following computations similar to those in \textbf{Step 1} of the first part of the proof, we obtain the existence of a constant $\mathsf{K}_{p,T,\mathfrak{C},\beta}$ such that
	$$
	\|(\bar{Y},\bar{Z},\bar{V},\bar{M})\|_{\mathcal{S}^{p,A}_\beta \times \mathcal{H}^{p}_\beta \times \mathfrak{L}^{p}_{\mu,\beta} \times \mathcal{M}^{p}_\beta}^p
	\leq \frac{\mathsf{K}_{p,T,\mathfrak{C},\beta}}{\varrho} \|(\bar{y},\bar{z},\bar{v},\bar{m})\|_{\mathcal{S}^{p,A}_\beta \times \mathcal{H}^{p}_\beta \times \mathfrak{L}^{p}_{\mu,\beta} \times \mathcal{M}^{p}_\beta}^p.
	$$
	It follows that for all $\varrho > \mathsf{K}_{p,T,\epsilon,\mathfrak{C},\beta}$, the mapping $\Psi$ is a strict contraction on $\mathcal{S}^{p,A}_\beta \times \mathcal{H}^p_\beta \times \mathfrak{L}^p_{\mu,\beta} \times \mathcal{M}^p_\beta$. Consequently, there exists a unique fixed point $(Y,Z,V)$ of $\Psi$ which, together with $K$, constitutes the unique $\mathbb{L}^p$-solution of the RBSDE \eqref{basic RBSDE} associated with the parameters $(\xi,f,L)$.
	
	To conclude the proof of Theorem \ref{thm2}, we now establish uniqueness. Let $(Y,Z,V,M,K)$ and $(Y',Z',V',M',K')$ be two $\mathbb{L}^p$-solutions of the RBSDE \eqref{basic RBSDE} associated with $(\xi,f,L)$ for some $p \in (1,2)$. Applying Itô's formula (Lemma \ref{Ito}) and using the Skorokhod condition, which implies $(Y_s - Y'_s)(dK_s - dK'_s) \leq 0$, we repeat the argument from \textbf{Step 1} of the first part of the proof. This yields the desired uniqueness.
\end{proof}

\begin{remark}
	The state process $(Y_t)_{t \leq T}$ of the RBSDE \eqref{basic RBSDE} can be represented as the solution to the following optimal stopping problem:
	$$
	Y_t = \esssup_{\tau \in \mathcal{T}_{[t,T]}} 
	\mathbb{E}\Bigg[ \int_{t}^{\tau} f(s,Y_s,Z_s,V_s)\,ds 
	+ L_\tau \mathds{1}_{\{\tau < T\}} 
	+ \xi \mathds{1}_{\{\tau = T\}} \;\Big|\; \mathcal{F}_t \Bigg], \qquad t \in [0,T].
	$$
\end{remark}

\end{document}